\documentclass[11pt,a4paper]{article} 

\usepackage{hyperref}
\hypersetup{colorlinks=true,linkcolor=red,citecolor=green,filecolor=magneta,urlcolor=cyan}

\usepackage{amsfonts,amssymb,amsmath,amsthm}
\usepackage[a4paper,hmargin={2.5cm,2.5cm},vmargin={2.3cm,2.3cm}]{geometry}

\usepackage{graphicx}

\usepackage[latin1]{inputenc}
\usepackage{color}
\newcommand{\Rset}{\mathbb{R}}
\newcommand{\Cset}{\mathbb{C}}

\newcommand{\Nset}{\mathbb{N}}

\newtheorem{thm}{Theorem}

\theoremstyle{definition}

\begin{document}

\title{Trapezoidal methods for fractional differential equations: theoretical and computational aspects 
\footnote{This is the post-print of the paper: R.Garrappa, ``Trapezoidal methods for fractional differential equations: theoretical and computational aspects'', \emph{Mathematics and Computers in Simulation} (Elsevier), April 2015, Volume 110, pp 96-112, doi: 10.1016/j.matcom.2013.09.012, avialable at \url{https://doi.org/10.1016/j.matcom.2013.09.012}. 
This work was supported by the GNCS-INdAM 2013 project ``Metodi numerici per equazioni differenziali alle derivate frazionarie''.}
} 

\author{
Roberto Garrappa \\
\small Universit\`a degli Studi di Bari, Dipartimento di Matematica, Bari, Italy\\ 
\small Member of the INdAM Research group GNCS \\
\small \texttt{roberto.garrappa@uniba.it}
}

\date{Submitted on March 12$^{\text{th}}$, 2013} 

\maketitle

\begin{abstract}
The paper describes different approaches to generalize the trapezoidal method to fractional differential equations. We analyze the main theoretical properties and we discuss computational aspects to implement efficient algorithms. Numerical experiments are provided to illustrate potential and limitations of the different methods under investigation.

\end{abstract}

%\maketitle


%\newproof{proof}{Proof}

\section{Introduction}

It is nowadays well established that several real--life phenomena are better described by fractional differential equations (FDEs), where the term \emph{fractional}, used for historical reasons, refers to derivative operators of any real positive order. Applications of FDEs are commonly found in bioengineering, chemistry, control theory, electronic circuit theory, mechanics, physics, seismology, signal processing and so on (e.g.,  \cite{BuenoOrovioKayGrauRodriguezBurrage2013,CaponettoMaionePisanoRapaiUsai2013,CafagnaGrassi2012,KilbasSrivastavaTrujillo2006,Garra2011,Magin2010,MachadoStefanescuTintareanuBaleanu2013}). We refer to \cite{MachadoKiryakovaMainardi2011} for an historical perspective on fractional calculus.

As a consequence of the growing interest for fractional order models, the analysis and the development of numerical methods for FDEs and partial differential equations with time or space fractional derivatives have become an active area of research (see, for instance, \cite{BurrageHaleKay2010,DiethelmFreed1999,DiethelmFordFreedLuchko2005,Diethelm2008,GarrappaPopolizio2011_JCAM,GarrappaPopolizio2013,KademLuchkoBaleanu2010,LiZeng2013,Lubich1986,MeerschaertTadjeran2004,MeerschaertTadjeran2006,MoretNovati2011,MoretPopolizio2012,Moret2013,YusteAcedo2005,YusteQuintanaMurillo2012}).

Some of the methods for FDEs are directly derived from methods for integral equations; this is the case, for instance, of product integration (PI) rules. A completely different approach, named as fractional linear multistep methods (FLMMs) \cite{Lubich1986}, has instead been specifically tailored to provide a solid theoretical background for the numerical treatment of FDEs.

Interestingly, these two approaches (PI and FLMM) lead to the same methods when the fractional order $\alpha$ converges to the nearest integer; indeed, different methods for FDEs are often generalizations of the same corresponding method for ordinary differential equations (ODEs).

Some important differences are however observed in the pure fractional case. For instance, the PI and FLMM generalizations of the implicit Euler method behave in a different way when applied to problems with discontinuous right--hand side \cite{DieciLopez2009,DieciLopez2011,DieciLopez2012} and only the method devised in the framework of FLMMs succeeds in stabilizing the solution on the switching surface, as expected from theoretical considerations, without generating unwanted spurious oscillations \cite{Garrappa2012_MCS}.

The presence of distinguishing characteristics in methods for FDEs developed from the same method for ODEs deserves to be more thoroughly investigated. Some differences are of theoretical interest (convergence, stability, etc.); nevertheless, the two approaches also present distinctive features of computational nature that affect the efficiency of the solution process.%with which the problem is solved.

The main aim of this paper is to describe and compare some methods for FDEs, highlight strengths and weaknesses and investigate the circumstances under which a certain method is preferable to some others. 

%Our analysis will be confined to those PI rules and FLMMs which can be interpreted, in some respects, as generalizations of the implicit second--order trapezoidal rule. This is a widely used method for ODEs which provides a satisfactory balance between accuracy, stability properties and computational complexity; indeed, it is an $A$--stable one--step method with the smallest possible error constant \cite{Lambert1991}. Our interest is in exploring whether these strengths are inherited also by the corresponding methods for FDEs.

The analysis will be confined to those PI rules and FLMMs that generalize, in some respects, the implicit second--order trapezoidal rule. This $A$--stable one--step method has the smallest possible error constant \cite{Lambert1991} and is appreciated since provides a satisfactory balance between accuracy, stability properties and computational complexity. Our interest is in exploring when these strengths are inherited by the corresponding methods for FDEs.

%The focusing of the attention on implicit methods is also motivated by their not frequent use in the context of FDEs: indeed, explicit schemes are usually preferred and, moreover, when implicit schemes are used, their use frequently occurs in the framework of predictor--corrector algorithms in which the implicit scheme is actually implemented as an explicit one, despite of the fact that stability reasons very often would require the use of fully implicit methods. 

There are also some more specific motivations for restricting the attention to second--order methods: a) PI rules very seldom converge with an order equal or greater than 2; b) the implementation of higher order FLMMs is usually computationally demanding; c) as for multistep methods for ODEs, order 2 represents a barrier for stability. 

Despite this apparently limited extent, we will however be able to include in our analysis a sufficiently large number of methods. 

The choice of focusing on implicit methods is also motivated by the fact that they have been so far very seldom investigated (very often, implicit methods are taken into account in the framework of predictor--corrector algorithms in which the implicit scheme is actually implemented as an explicit one) and explicit methods are unduly preferred in most of the works on FDEs, regardless of stability issues which discourage their use.

%Focusing the attention on implicit methods is also useful since explicit methods are unduly preferred (regardless of stability issues that would suggest otherwise) in most of the works on FDEs and very few works are devoted to implicit schemes (very often, implicit methods are taken into account in the framework of predictor--corrector algorithms in which the implicit scheme is actually implemented as an explicit one). 

%The focusing of the attention on implicit methods finds its reasons on the lack of .. not frequent use with FDEs; not only explicit methods are in most cases preferred (despite stability issues that would suggest otherwise) but the use of implicit schemes occurs nearly always in the framework of predictor--corrector algorithms in which the implicit scheme is actually implemented as an explicit one. 

%The restriction to second--order methods is more involved and concerns with different aspects: PI rules, as we will see later, very seldom converge with an order equal or greater than 2; the implementation of higher order FLMMs can prove to be highly computationally expensive; as in the case of ordinary differential equations (ODEs), order 2 of convergence is a barrier for stability. Despite this apparently limited extent, we will however be able to include in our analysis a sufficiently large number of methods. 

The paper is organized as follows. In Section \ref{S:FDE} we present some basic definitions and properties concerning FDEs. Section \ref{S:ProductIntegration} is devoted to describe PI rules and discuss different implementations on uniform and graded grids. In Section \ref{S:FLMM} we present FLMMs and we describe three different generalizations of the implicit trapezoidal rule; the use of an effective formula for the evaluation of the weights is also proposed. In Section \ref{S:LinearStability} linear stability is investigated and some implementation details are discussed in Section \ref{S:ImpementationIssues}. Finally, in Section \ref{S:NumericalExperiments} the results of some numerical tests are presented to compare the methods under investigation and some concluding remarks are presented at the end of the paper.

\section{Fractional differential equations}\label{S:FDE}
Let us consider the initial value problem for a system of FDEs in the form
\begin{equation}\label{eq:FDE}
	\left\{ \begin{array}{l}
		{}^C D^{\alpha}_{t_{0}} y(t) = f(t,y(t)) \\
		y^{(k)}(t_{0}) = y_{0,k} , \quad
		k = 0, \dots m -1
	\end{array} \right. ,
\end{equation}
where $m=\left\lceil \alpha \right\rceil$ is the smallest integer such that $m>\alpha$,
$f:[t_{0},T] \times \Rset^{q} \to \Rset^{q}$ is assumed to be sufficiently smooth, $y : [t_{0},T] \to \Rset^{q}$ is the unknown solution and $y^{(k)}$ denotes classical derivative of integer order $k$. 

The fractional derivative operator ${}^C D^{\alpha}_{t_{0}}$ is introduced in this paper according to the Caputo's definition  \cite{Diethelm2010,Mainardi2010,Podlubny1999}
\[
	{}^C D^{\alpha}_{t_{0}} y(t) = I_{t_{0}}^{m - \alpha} y^{(m)}(t) ,
%	= \frac{1}{\Gamma(m-\alpha)} \int_{t_{0}}^{t} (t-s)^{m-\alpha-1} y^{(m)}(s) ds
\]
with $I_{t_{0}}^{\beta}$ the \emph{Riemann--Liouville} (RL) integral of order $\beta$ on the interval $[t_0, t]$ 
\begin{equation}\label{eq:RiemannLiouvilleIntegral}
	I_{t_{0}}^{\beta} g(t) = \frac{1}{\Gamma(\beta)} \int_{t_{0}}^{t} (t-s)^{\beta-1} g(s) ds .
\end{equation}

A classical result in fractional calculus \cite{Diethelm2010} allows to compute, for any $\alpha \in \Rset$, the RL integral of the power function as
\begin{equation}\label{eq:RL_Integral_Power}
	I_{t_{0}}^{\alpha} (t-t_{0})^{\nu} = \frac{\Gamma(\nu+1)}{\Gamma(\alpha+\nu+1)} (t-t_{0})^{\alpha+\nu}
	, \quad
	\nu > -1
	, \quad
	t \ge t_{0} .
\end{equation}

It is a well known result that (\ref{eq:FDE}) can be equivalently written in terms of the integral equation 
\begin{equation}\label{eq:FDE_RL_1}
	y(t) = T_{m-1}(t) + I_{t_{0}}^{\alpha} f(t,y(t)) ,
\end{equation}
where $T_{m-1}(t)$ is the Taylor expansion of $y(t)$ centered at $t_{0}$ 
\[
	T_{m-1}(t) = \sum_{k=0}^{m -1} \frac{(t-t_{0})^{k}}{k!} y_{0}^{(k)} .
\]

One of the major difficulties in the treatment of FDEs derives from the lack of smoothness of the true solution on the whole interval of integration. Indeed, as proved in \cite{Lubich1983}, the true solution of (\ref{eq:FDE}) possesses an expansion in mixed (integer and real) powers of the form
\begin{equation}\label{eq:FDE_Expansion_True}
	y(t) = T_{m-1}(t) +  \! \! \sum_{\nu \in {\cal A}_{p,m}} \! \! (t-t_{0})^{\nu} Y_{\nu}  + {\cal O} \bigl( (t-t_{0})^{p} \bigr) 
	, \quad t \to t_{0}, 
\end{equation}
with $p > m$, ${\cal A}_{p,m} = {\cal A}_{p} - \bigl\{ 0, \dots , m-1 \bigr\}$ and 
\[
	{\cal A}_{p} = \bigl\{\nu \in \Rset \, | \, \nu = i+j\alpha, \, i,j\in \Nset, \, \nu < p \bigr\} ,
\]
with the consequence that the derivatives of $y(t)$ are unbounded at $t_{0}$.

\section{Product integration rules}\label{S:ProductIntegration}

The main idea behind PI rules \cite{Young1954} is to evaluate the integral in (\ref{eq:FDE_RL_1}) by approximating the vector field $f$ with suitable polynomials. For the integer--order case this is the way in which Adams multistep methods are devised. 

To generalize the $1$--step trapezoidal rule consider a (not necessarily uniform) grid $\bigl\{t_{0},t_{1},\dots,t_{N}\bigr\}$ on the whole interval of integration $[t_{0},T]$ and, in the piecewise decomposition of (\ref{eq:FDE_RL_1})
\[
	y(t) = T_{m-1}(t) + \frac{1}{\Gamma(\alpha)} \sum_{j=0}^{n-1} \int_{t_{j}}^{t_{j+1}} (t-s)^{\alpha-1} f(s,y(s)) ds
	, \quad
	t \ge t_{n} ,
\]
replace $f$ in each subinterval by the first--degree polynomial interpolant
\[
	p_{j}(s) = f_{j+1} + \frac{s-t_{j+1}}{h_{j}} \bigl( f_{j+1} - f_{j} \bigr) 
	, \quad  s \in [t_{j},t_{j+1}],
\]
where $h_{j} =  t_{j+1}-t_{j}$, $f_{j}=f(t_{j},y_{j})$ and $y_{j}$ is the numerical approximation to the solution $y(t_{j})$. After exactly solving the integrals
\begin{eqnarray*}
	\frac{1}{\Gamma(\alpha)} \int_{t_{j}}^{t_{j+1}} (t_{n}-s)^{\alpha-1} ds
	&=& I_{n,j}^{(0)} - I_{n,j+1}^{(0)} ,\\
	\frac{1}{\Gamma(\alpha)} \int_{t_{j}}^{t_{j+1}} (t_{n}-s)^{\alpha-1} \frac{(s - t_{j+1})}{h_{j}} ds
	&=& \frac{I_{n,j}^{(1)}}{h_{j}} - I_{n,j}^{(0)} - \frac{I_{n,j+1}^{(1)}}{h_{j}} , \\
\end{eqnarray*}
where, for shortness, we denoted 
\[
	I_{n,j}^{(k)} = \frac{1}{\Gamma(\alpha)} \int_{t_{j}}^{t_{n}} (t_{n}-s)^{\alpha-1} (s-t_{j})^k ds = 
	\frac{(t_{n}-t_{j})^{\alpha+k}}{\Gamma(\alpha+k+1)} ,
\]
it is elementary, after a simple change of the summation index, to obtain the numerical approximation of the solution $y(t)$ at $t=t_{n}$ as 
\begin{equation}\label{eq:PI_General}
	y_{n} = T_{m-1}(t_{n}) 
	+ w_{n} f_{0} 
	+ \sum_{j=1}^{n} b_{n,j} f_{j} ,
\end{equation}
where
\[
	\left\{\begin{array}{l}
		\displaystyle
		w_{n} = I_{n,0}^{(0)} - \frac{I_{n,0}^{(1)}}{h_{0}} + \frac{I_{n,1}^{(1)}}{h_{0}} , \\
		\displaystyle
		b_{n,j} = \frac{I_{n,j-1}^{(1)} - I_{n,j}^{(1)}}{h_{j-1}} - \frac{I_{n,j}^{(1)}- I_{n,j+1}^{(1)}}{h_{j}} 
		, \quad
		j = 1,\dots,n-1
		, \quad
		b_{n,n} =  \frac{I_{n,n-1}^{(1)}}{h_{n-1}} .\\
	\end{array} \right.
\]

According to the analysis presented in \cite{CameronMcKee1985}, this scheme is convergent of second order whenever the true solution is assumed sufficiently smooth. Unfortunately, this is a rather optimistic result compared to the expansion (\ref{eq:FDE_Expansion_True}) of $y(t)$: due to the presence of unbounded derivatives at $t_{0}$ the convergence rate of (\ref{eq:PI_General}) is actually lower than the theoretical value $2$, as it can be formally stated by a direct application of the main result by Dixon in \cite{Dixon1985}.

\begin{thm}\label{thm:PI_Convergence}
Let $f$ be Lipschitz continuous with respect to the second variable and $y_{n}$ the numerical approximation obtained by applying the PI trapezoidal rule (\ref{eq:PI_General}) on the interval $[t_{0},T]$. There exists a constant $C=C(T-t_{0})$, which does not depend on $h$, such that
\[
	\left \| y(t_{n}) - y_{n} \right \| \le C \bigl( t_{n}^{\alpha-1} h^{1+\alpha} + h^{2} \bigr) 
	, \quad
	h = \max_{j=0,\dots,n-1} h_{j}.
\]
\end{thm}

Unless very restrictive (and usually unrealistic) assumptions are made to ensure the smoothness of $f(t,y(t))$ also at the left endpoint of $[t_{0},T]$, the rate by which the error reduces to $0$ is therefore, away from $t_{0}$, proportional to $h^{1+\alpha}$ when $0<\alpha<1$. The same drop in the convergence order occurs also when higher degree polynomials are employed; for this reason PI rules of higher order are never taken into account for FDEs when $0<\alpha<1$. When $\alpha>1$ order $2$ of convergence is instead obtained.

%The problem of finding suitable conditions on the vector field in order to ... .... ...  An attempt to find conditions of the vector field such that higher order can be obtained has been investigated in ... ... Champing Li o forse Deng Wehuwhua ... ... verificare ....  ... ... ...  Vedere anche Weihua Deng, Smoothness and stability of the solutions for nonlinear fractional differential equations, Nonlinear Anal., 72(2010) 17681777

\subsection{Uniform meshes}\label{SS:PI_Uniform}

Using (\ref{eq:PI_General}) on uniform grids $t_{n}=t_{0}+nh$, with $h>0$ a constant step--size, offers some advantages of computational type. The evaluation of the weights is indeed less expensive and just involves the computation of real powers of integer numbers. Furthermore, (\ref{eq:PI_General}) presents a convolution structure, in which each weight $b_{n,j}$ depends just on $n-j$, and it can be simplified as 
\begin{equation}\label{eq:PI_General_Uniform}
	y(t_{n}) = T_{m-1}(t_{n}) 
	+ \frac{h^{\alpha}}{\Gamma(\alpha+2)} \left( \tilde{w}_{n} f_{0} 
	+ \sum_{j=1}^{n} \tilde{b}_{n-j} f_{j} \right) ,
\end{equation}
where now 
\[
	\left\{\begin{array}{l}
		\tilde{w}_{n} = (\alpha+1-n)n^{\alpha} + (n-1)^{\alpha+1} \\
		\tilde{b}_{0} =  1 , \quad
		\tilde{b}_{n} = (n-1)^{\alpha+1} - 2n^{\alpha+1} + (n+1)^{\alpha+1}
		, \quad n=1,2,\dots \\
	\end{array}\right. .
\]

Formulation (\ref{eq:PI_General_Uniform}) of the trapezoidal PI rule is well--known since it is a building component of a widely used predictor--corrector algorithm \cite{DiethelmFreed1999}.% for which an optimized Matlab implementation is available \cite{Garrappa_FDE12}. 

The main advantage of the convolution structure is however that (\ref{eq:PI_General_Uniform}) can be implemented in a fast way by means of FFT algorithms \cite{HairerLubichSchlichte1985}.

\subsection{Graded meshes}\label{SS:PI_Graded}

Unlike FLMMs, PI rules can be implemented on non uniform grids, with the aim of overcoming the order reduction stated in Theorem \ref{thm:PI_Convergence}. This issue has been investigated, in the context of weakly singular integral equations, by several authors \cite{Brunner1985,HoogWeiss1974} and graded grids are highly praised to this purpose.

%By clustering the grid--points at the left boundary of the interval of integration, graded grids are able to cope with the lack of regularity and counteract the inefficiency of polynomials to satisfactorily approximate real powers. 

By clustering the grid--points at the left boundary of the integration interval, graded grids are able to cope with the lack of regularity and counteract the inefficiency of polynomials to satisfactorily approximate real powers. 

A graded mesh on the interval $[t_{0},T]$ is defined according to
\[
	t_{n} = t_{0} + \left( \frac{n}{N} \right)^{r} (T-t_{0}) 
	, \quad
	n = 0, 1, \dots, N,
\]
where $r>1$ is a suitable grading exponent. After denoting $h_{0} = t_{1} - t_{0} = (T-t_{0})/{N^r}$ and observing that $h_{j} = t_{j+1} - t_{j} = h_{0} \bigl( (j+1)^r - j^r\bigr)$, a simple manipulation allows to reformulate (\ref{eq:PI_General}) according to 
\[
	y(t_{n}) = T_{m-1}(t_{n}) 
	+ \frac{h_{0}^{\alpha}}{\Gamma(\alpha+2)} \left( \hat{w}_{n} f_{0} 
	+ \sum_{j=1}^{n} \hat{b}_{n,j} f_{j} \right) ,
\]
where
\[
	\left\{ \begin{array}{l}
		\displaystyle
		\hat{w}_{n} = (n^r -1)^{\alpha+1} - n^{r \alpha} (n^{r} - \alpha -1) , \\
		\displaystyle
		\hat{b}_{n,j} = \phi_{n,j}(\alpha,r) - \phi_{n,j+1}(\alpha,r) 
		, \quad
		j = 1,\dots,n-1
		, \quad
		\hat{b}_{n,n} =  (n^{r} - (n-1)^{r})^{\alpha} \\
	\end{array}\right.
\]
and
\[
	\phi_{n,j}(\alpha,r) = \frac{(n^r-(j-1)^r)^{(\alpha+1)} - (n^r-j^r)^{(\alpha+1)}}{j^r-(j-1)^r}	.
\]

The following result \cite{Brunner1985,Brunner2004} on the second order convergence of graded grids for linear FDEs can be readily extended to more general equations.

\begin{thm}
Let $f(t,y(t)) = K(t)y(t)$, $K \in {\cal C}^{2}([t_{0},T])$. For $r=2/\alpha$ it is 
\[
	\left \| y(t_{n}) - y_{n} \right \| \le C_{\alpha} h^{2} 
	, \quad
	h = \! \! \! \! \max_{j=0,\dots,n-1} \! \! \! \! h_{j} = h_{n-1}.
\]
\end{thm}

The use of PI rules on nonuniform grids is more demanding; not only the weights must be recomputed at each step but, since the absence of a convolution structure, it is not possible to employ fast algorithms. By means of some numerical experiments we will verify, at the end of the paper, whether or not this further amount of computation is, in some way, balanced by the enhancement in the accuracy.

%Our goal is therefore to verify, by means of numerical experiments at the end of the paper, whether adding this further amount of computation is, in some way, balanced by the enhancement in the accuracy.

\section{Fractional linear multistep methods }\label{S:FLMM}

In the pioneering work on discretized fractional calculus \cite{Lubich1986}, Lubich proposed an elegant and effective strategy for classical linear multistep methods (LMMs) originally devised for integer--order ODEs.  

Although the potentials of FLMMs have been recognized in several papers (see, for instance, the review article by Gorenflo \cite{Gorenflo1997}), their use for practical computation has not captured so far widespread attention, perhaps due to some difficulties in explicitly formulating FLMMs.

The key aspect in FLMMs is the approximation of the RL integral (\ref{eq:RiemannLiouvilleIntegral}) by means of the convolution quadrature 
\begin{equation}\label{eq:ConvolutionQuadrature}
	I_{h}^{\beta} g(t_{n}) 
	= h^{\beta} \sum_{j=0}^{n} \omega_{n-j} g(t_{j}) + h^{\beta} \sum_{j=0}^{s} w_{n,j} g(t_{j})
\end{equation}
on uniform grids $t_{n}=t_{0}+nh$, $h>0$, and where convolution and starting quadrature weights $\omega_{n}$ and $w_{n,j}$ do not depend on $h$. 

Starting weights $w_{n,j}$ play an important role, especially in the first part of the integration interval, in order to cope with the possibly singular character of the integrand function at $t_{0}$; we will discuss the problem of their evaluation later on.

Convolution quadrature weights $\omega_{n}$ are the main components of the quadrature rule and characterize the specific FLMM. They are obtained starting from any LMM for ODEs
\[
	\sum_{j=0}^{k} \rho_{j} y_{n-j} = \sum_{j=0}^{k} \sigma_{j} f(t_{n-j},y_{n-j}) ,
\]
being $\rho(z)=\rho_{0}z^{k} + \rho_{1}z^{k-1} + \dots + \rho_{k}$ and $\sigma(z)=\sigma_{0}z^{k} + \sigma_{1}z^{k-1} + \dots + \sigma_{k}$ the characteristic polynomials. As shown in \cite{Wolkenfelt1979}, LMMs can be equivalently reformulated as (\ref{eq:ConvolutionQuadrature})  (but with $\beta=1$ for ODEs) with weights $\omega_{n}$ obtained as the coefficients of the formal power series (FPS)
\[
	\omega(\xi)  = \sum_{n=0}^{\infty} \omega_{n} \xi^{n}
	, \quad
	\omega(\xi) = \frac{\sigma(1/\xi)}{\rho(1/\xi)} 
\]
and the function $\omega(\xi)$ goes therefore under the name of \emph{generating function} of the LMM. The idea proposed in \cite{Lubich1986} is to generate a quadrature rule for fractional problems (\ref{eq:FDE}) by evaluating the convolution weights as the coefficients in the FPS of the fractional--order power of the generating function
\begin{equation}\label{eq:PowerGeneratingFunction}
	\omega_{\beta}(\xi) = \sum_{n=0}^{\infty} \omega_{n} \xi^{n}
	, \quad \omega_{\beta}(\xi) = \left( \frac{\sigma(1/\xi)}{\rho(1/\xi)} \right)^{\beta} .
\end{equation}

Methods of this kind, named as FLMMs, when applied to (\ref{eq:FDE}) read as
\begin{equation}\label{eq:FLMM}
	y_{n} = T_{m-1}(t_{n}) + h^{\alpha} \sum_{j=0}^{s} w_{n,j} f_{j} + h^{\alpha} \sum_{j=0}^{n} \omega_{n-j} f_{j} 
\end{equation}
and their convergence properties are stated in the following result \cite{Lubich1985,Lubich1986}.

\begin{thm}\label{thm:FLMM_Convergence}
Let $(\rho,\sigma)$ be a stable and consistent of order $p$ implicit LMM with the zeros of $\sigma(\xi)$ having absolute value $\le1$. The FLMM (\ref{eq:FLMM}) is convergent of order $p$.
\end{thm}

One of the major difficulties in FLMMs (\ref{eq:FLMM}) is in evaluating the weights $\omega_{n}$ as the coefficients in the FPS (\ref{eq:PowerGeneratingFunction}). Although some sophisticated algorithms are available for manipulating FPS \cite{Henrici1979}, for most of the methods an efficient tool is the J.C.P. Miller formula stated by the following theorem \cite[Theorem 1.6c]{Henrici1974}.

\begin{thm}\label{thm:PowerPowerSeries}
Let $\varphi(\xi) = 1 + \sum_{n=1}^\infty a_n \xi^n$ be a FPS. Then for any $\beta \in \Cset$, 
\[
	\bigl( \varphi(\xi) \bigr)^{\beta} = \sum_{n=0}^\infty v_n^{(\beta)} \xi^n ,
\]
where coefficients $v_n^{(\beta)}$ can be recursively evaluated as
\[
	v_0^{(\beta)}=1 , \quad 
	v_n^{(\beta)} = \sum_{j=1}^n \left( \frac{(\beta+1)j}{n} - 1 \right) a_j v^{(\beta)}_{n-j} .
\]
\end{thm}

The Miller formula allows to evaluate, in the most general case, the first $N$ coefficients of $\bigl( \varphi(\xi) \bigr)^\beta$ with a number of operations proportional to $N^{2}$. In most cases of practical interest, this computational effort is actually reduced by an order of magnitude. This is the case of the evaluation of the coefficients $\bigl\{ \omega_{n}^{(\beta)} \bigr\}_{n\in \Nset}$ in the FPS of $(1 \pm \xi)^{\beta}$; indeed, since $a_{1}=\pm 1$ and $a_{2}=a_{3}=\dots=0$, it is elementary to see that the application of Theorem \ref{thm:PowerPowerSeries} leads to
\begin{equation}\label{eq:FPS_Recursive}
	\omega_{0}^{(\beta)} = 1, \quad \omega_{n}^{(\beta)} = \pm \left( \frac{\beta+1}{n} - 1 \right) \omega_{n-1} ,
\end{equation}
involving just $2N$ multiplications and $N$ additions.

\subsection{Fractional trapezoidal rule}\label{SS:FracTrapRule}
In order to generalize to FDEs the classical trapezoidal rule, let us first consider its classical formulation for ODEs
\begin{equation}\label{eq:TrapezoidalODE}
	y_{n+1} - y_{n} = \frac{h}{2} \bigl( f_{n} + f_{n+1} \bigr) ,
\end{equation}
with characteristic polynomials $\rho(z) = z-1$ and $\sigma(z) = (z+1)/2$ and generating function 
\[
	\omega(\xi) = \frac{\sigma(1/\xi)}{\rho(1/\xi)}
	= \frac{(1+\xi)}{2(1-\xi)} 
	= \frac{1}{2} \left( 1 + 2 \sum_{n=1}^{\infty} \xi^n \right) .
\]

To practically compute the first $N$ weights in the FPS of $\bigl( \omega(\xi)\bigr)^{\alpha}$ we advise against the direct application of the Miller formula to the infinite series in $\omega(\xi)$ since it would involve a number of operations proportional to $N^{2}$. It is instead preferable to first apply twice Theorem \ref{thm:PowerPowerSeries} to $\bigl(1+\xi\bigr)^{\alpha}$ and $\bigl(1-\xi\bigr)^{-\alpha}$, by means of  (\ref{eq:FPS_Recursive}), and hence evaluate the coefficient of their product by an FFT algorithm, with a number of operations proportional to $3N \log_{2} 4 N$ when $N$ is a power of $2$ \cite{Henrici1979}. This last task can be efficiently performed by means of the very few lines of Matlab code
\begin{verbatim}
x = fft([omega1,zeros(size(omega1))]) ;
y = fft([omega2,zeros(size(omega2))]) ;
omega = ifft(x.*y) ; omega = omega(1:N) ;
\end{verbatim}
where {\tt omega1} and {\tt omega2} are the row arrays with coefficients of $\bigl(1+\xi\bigr)^{\alpha}$ and $2^{-\alpha}\bigl(1-\xi\bigr)^{-\alpha}$ respectively, evaluated in a fast way thanks to (\ref{eq:FPS_Recursive}).

\subsection{Newton--Gregory formula}\label{SS:NewtonGregoryFormula}

The trapezoidal rule (\ref{eq:TrapezoidalODE}) belongs to the more general family of $k$--step Adams--Moulton methods
\[
	y_{n+1} - y_{n} = h \sum_{i=0}^{k} \gamma_{i} \nabla^{i} f_{n+1} ,
\]
with $\nabla^{i} f_{n+1}$ the classical backward differences (i.e., $\nabla^{0} f_{n+1} = f_{n+1}$ and $\nabla^{i+1} f_{n+1} = \nabla^{i} f_{n+1} - \nabla^{i} f_{n}$); the coefficients $\gamma_{i}$ are the first $k+1$ terms in the truncation 
\[
	\bigl[ G(\xi) \bigr]_{k} = \gamma_{0} + \gamma_{1} (1-\xi) + \dots + \gamma_{k} (1-\xi)^k 
\]
of the power series expansion, with respect to $1-\xi$, of $G(\xi) = (\xi-1)/\ln \xi$ and the resulting generating function is $\omega(\xi) = \bigl[ G(\xi) \bigr]_{k}/(1-\xi)$. 

It is easy to see that $\gamma_{0}=1$ and $\gamma_{1} = -\frac{1}{2}$, thus immediately leading to (\ref{eq:TrapezoidalODE}) when $k=1$. Therefore, a straightforward generalization to fractional--order problems by means of the $\alpha$--power of $\omega(\xi)$ leads, in the case of the $1$--step method, to the trapezoidal rule discussed in the previous subsection.

It has been shown in \cite{Lubich1983} that alternative methods, named as \emph{Newton--Gregory formulas}, can be devised by first evaluating the $\alpha$ power of $G(\xi)$ and successively truncating to the first $k+1$ terms as
\[
	\omega_{\alpha}(\xi) = \frac{\bigl[ \left( G(\xi) \right)^{\alpha} \bigr]_{k}}{(1 - \xi)^{\alpha}} 	.
\]

The interchange between the truncation series symbol $[ \, \cdot \, ]_{k}$ and the $\alpha$--power is important and generates a different family of methods. By operating for simplicity the change of variable $\tau = 1-\xi$, 
\[
	\left(G(\xi) \right)^{\alpha} = \left(G(1-\tau) \right)^{\alpha} 
	%= \left( \frac{-z}{\ln(1-z)} \right)^{\alpha}
	= \left( -\frac{\ln(1-\tau)}{\tau} \right)^{-\alpha}
	= \left( 1 + \sum_{n=1} \frac{\tau^{n}}{n+1}  \right)^{-\alpha}
	= \sum_{n=0}^{\infty} \bar{\gamma}^{\alpha}_{n} \tau^{n} .
\]

The first coefficients $\bar{\gamma}^{\alpha}_{n}$ can be easily evaluated thanks again to Theorem \ref{thm:PowerPowerSeries} (see \cite{GaleoneGarrappa2008} for the explicit formulation of the first coefficients).

Since $\bar{\gamma}^{\alpha}_{0} = 1$ and $\bar{\gamma}^{\alpha}_{1} = - \frac{\alpha}{2}$, the generating function of the FLLM obtained by the $1$--step Adams method is therefore
\[
	\omega_{\alpha}(\xi) =  \frac{1 - \frac{\alpha}{2}(1-\xi) }{(1 - \xi)^{\alpha}} 
	= (1 - \xi)^{-\alpha} \left( 1 - \frac{\alpha}{2}(1-\xi) \right) 
\]
and, after a simple manipulation, it is possible to evaluate the first $N$ coefficients in $\omega_{\alpha}(\xi)$, with a cost proportional to $N$, as
\[
	\omega_{0} = 1-\frac{\alpha}{2}
	, \quad
	\omega_{n} = \left(1-\frac{\alpha}{2}\right) \omega_{n}^{(-\alpha)} + \frac{\alpha}{2} \omega_{n-1}^{(-\alpha)} ,
\]
with $\omega_{n}^{(-\alpha)}$ the coefficients in the FPS $(1 - \xi)^{-\alpha}$ evaluated by means of (\ref{eq:FPS_Recursive}).

\subsection{Fractional BDF formula}\label{SS:FracBDF}

The second order BDF formula for ODEs is given by
\begin{equation}\label{eq:BDF2}
	y_{n+2} - \frac{4}{3} y_{n+1} + \frac{1}{3} y_{n} = \frac{2h}{3} f_{n+2} ,
\end{equation}
with characteristic polynomials $\rho(z)= (z^2 - 4/3 z + 1/3)$ and $\sigma(z)=2/3z^2$ and generating function 
\[
	\omega(\xi) = \frac{\sigma(1/\xi)}{\rho(1/\xi)}
%	= \frac{1}{(3/2-2\xi+\xi^2/2)} = \frac{2}{(3-\xi)(1-\xi)}
	= \frac{2}{3 \bigl( 1 - 4\xi/3 + \xi^2/3 \bigr)} .
\]

BDF method (\ref{eq:BDF2}) can be considered in no way as coming from the trapezoidal rule; it is indeed obtained by differentiating a second--degree interpolant polynomial. However, it is equally of interest in our analysis due to its second order of convergence as the trapezoidal rule and the excellent stability properties that we are interested in verifying also for FDEs.

%Although (\ref{eq:BDF2}) can not be considered as originated from the trapezoidal rule, it is the same of interest in our analysis since its second order of convergence.

The coefficients in the FPS of $\bigl( \omega(\xi) \bigr)^{\alpha}$ can be profitably evaluated again by the Miller formula applied to $1-4\xi/3 + \xi^2/3$, with negative power $-\alpha$, to first produce
\[
	\tilde{\omega}_{0} = 1
	, \quad
	\tilde{\omega}_{1} = \frac{4}{3} \alpha \tilde{\omega}_{0}
	, \quad
	\tilde{\omega}_{n} = \frac{4}{3} \left( 1 + \frac{\alpha-1}{n} \right) \tilde{\omega}_{n-1} 
	+ \frac{4}{3} \left( \frac{2(1-\alpha)}{n} - 1 \right) \tilde{\omega}_{n-2} 
\]
and hence by computing $\omega_{n} = 2^{\alpha}\tilde{\omega}_{n}/3^{\alpha}$ with an overall number of operations that still remains proportional to $N$.

\section{Linear stability analysis}\label{S:LinearStability}

The analysis of stability is of primary importance to understand the qualitative behavior of a numerical method. By assuming the existence of a stable steady state, stability analysis aims to verify whether or not approximated trajectories move towards the steady state.

As usual, stability issues for FDEs can be addressed by studying the scalar linear test equation
\begin{equation}\label{eq:FDE_Linear}
	{}^C D^{\alpha}_{t_{0}} y(t) = \lambda y(t)
	, \quad \lambda \in \Cset ,
\end{equation}
whose steady--state is at $y = 0$. The generalization of a well--known result in the theory of ODEs can be done by means of the following result \cite{Lubich1985,Matignon1998}. 

\begin{thm}\label{thm:FDEStability}
Let $\alpha>0$. The steady--state $y=0$ of (\ref{eq:FDE_Linear}) is stable if and only if $\lambda \in \Sigma_{\alpha}$, where $\Sigma_{\alpha} =  \left\{ s \in \Cset \, \, : \, \, \bigl|\arg(s)\bigr| > {\alpha \pi}/{2} \right\}$.
\end{thm}

\begin{figure}[htb]
	\centering
	\begin{tabular}{c@{\hspace{0.8cm}}c}
		\includegraphics[width=0.38\textwidth,height=0.40\textwidth]{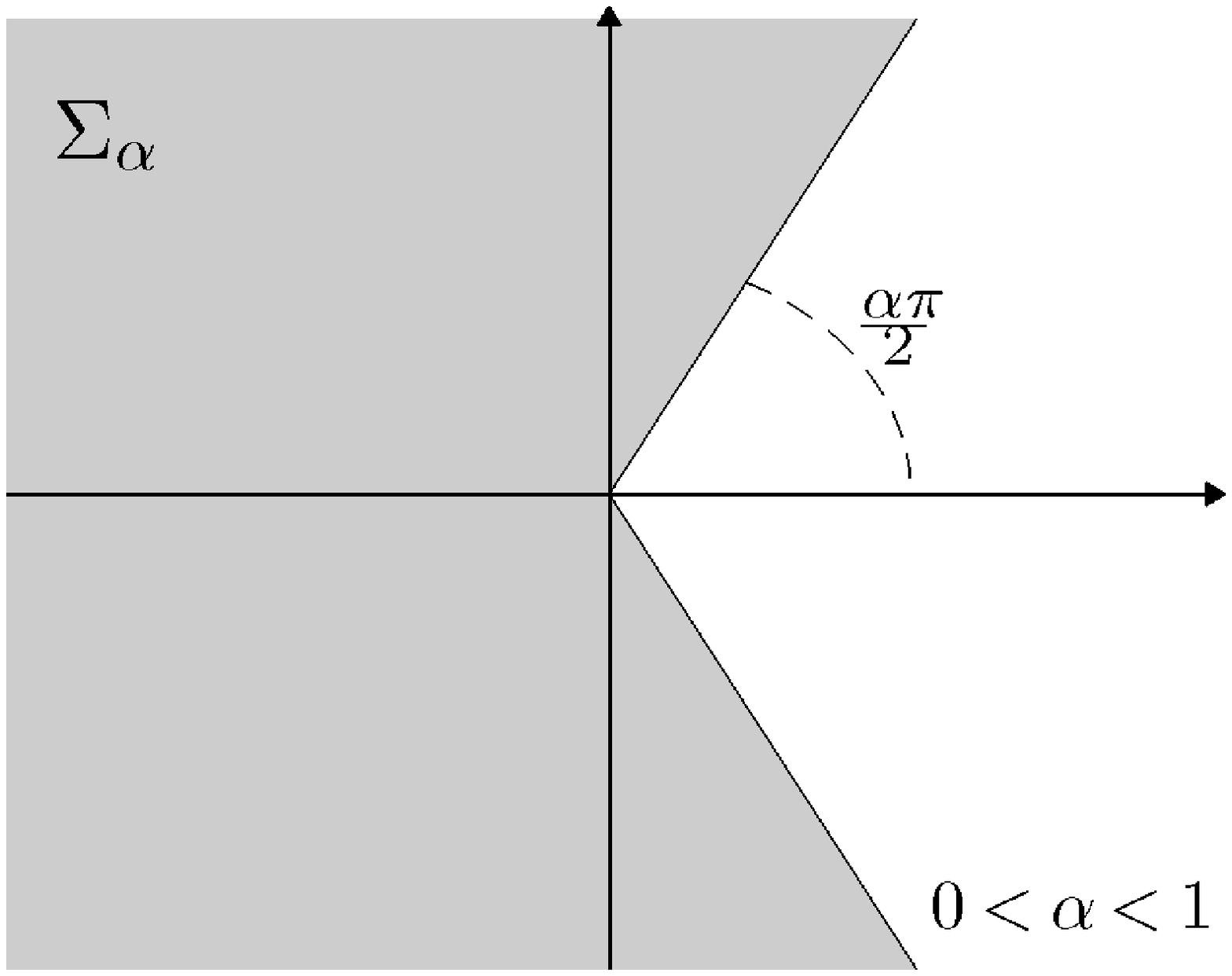} 
		&
		\includegraphics[width=0.38\textwidth,height=0.40\textwidth]{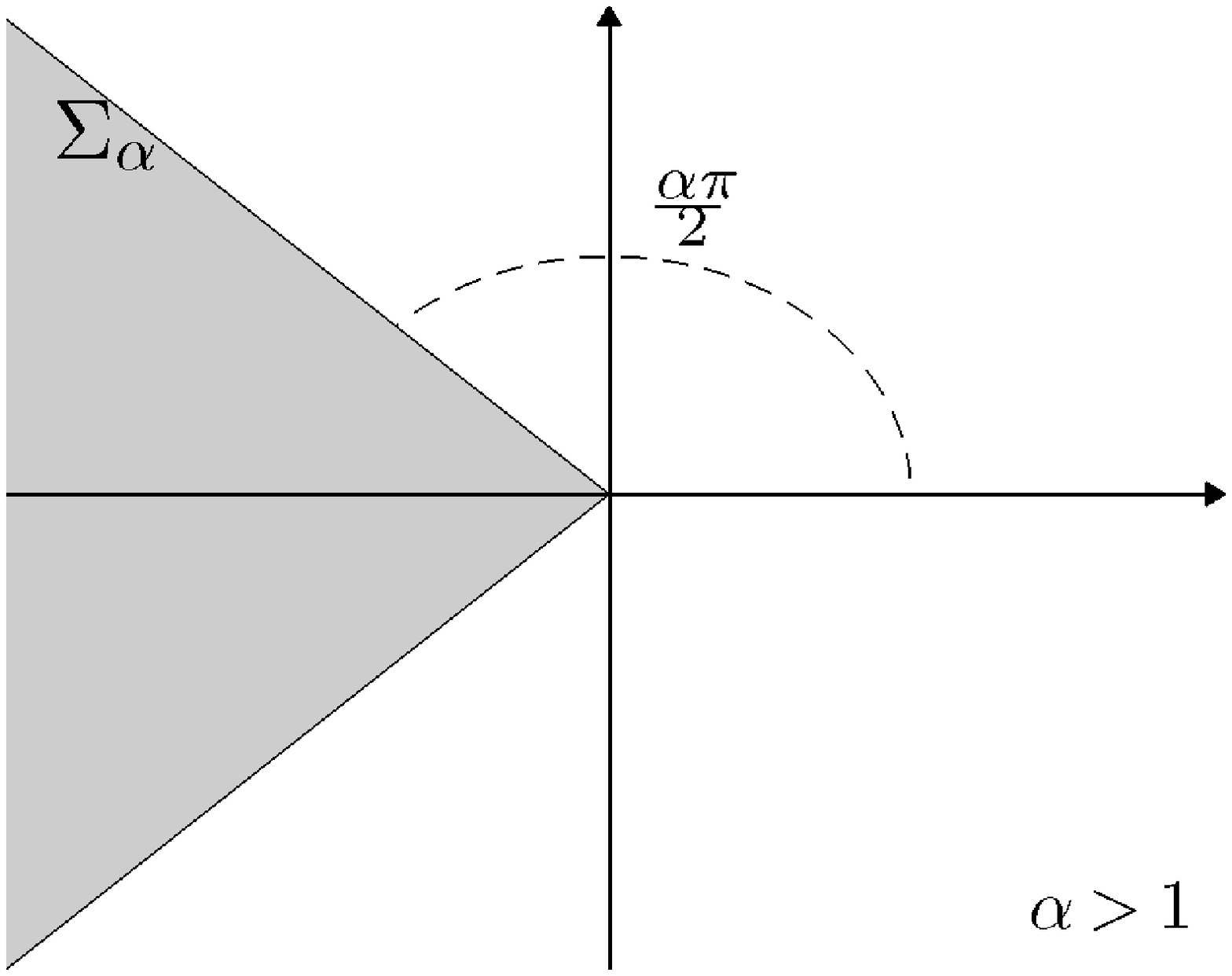} 
	\end{tabular}
\caption{Stability sectors $\Sigma_{\alpha}$ for $0<\alpha<1$ (left) and $\alpha>1$ (right)}
\label{fig:Fig_StabilitySector}
\end{figure}

In Figure \ref{fig:Fig_StabilitySector} we present the sector of stability $\Sigma_{\alpha}$ for two exemplifying cases. It is immediate to verify that $\Sigma_{\alpha}$ coincides with the left complex semi--plane $\Cset^{-}$ when $\alpha=1$, thus generalizing the classical result for ODEs.

%In Figure \ref{fig:Fig_StabilitySector} we present the sector of stability $\Sigma_{\alpha}$ for two exemplifying cases. It is an immediate task to verify that, whenever $\alpha=1$, $\Sigma_{\alpha}$ coincides with the left complex semi--plane $\Cset^{-}$ representing the stability region of classical ODEs.

%In the case of a linear systems of ODEs, the same requirments are ... ... for the eigenvalue of the matrix of the ... system.

When PI rules or FLMMs are applied to (\ref{eq:FDE_Linear}), the numerical solution $\bigl\{y_{n}\bigr\}_{n\in\Nset}$ is evaluated by means of a convolution quadrature
\begin{equation}\label{eq:GeneralConvolutionQuadrature}
	y_{n} = g_{n} + z_{\alpha} \sum_{j=0}^{n} \omega_{n-j} y_{j} 
	, \quad z_{\alpha} = h^{\alpha} \lambda .
\end{equation}

The problem of determining the stability region $S_{\alpha} \subseteq \Cset$ such that $y_{n} \to 0$ whenever $z_{\alpha} \in S_{\alpha}$ has been discussed in \cite{Lubich1985,Lubich1986_IMA} and the main result can be summarized in the following theorem.

\begin{thm}\label{thm:FED_StabilityRegion}
Let $\alpha>0$. Under the assumption that the sequence $\bigl\{g_{n}\bigr\}_{n}$ is convergent and that the quadrature weights $\omega_{n}$ satisfy
\[
	\omega_{n} = \frac{n^{\alpha-1}}{\Gamma(\alpha+1)} + u_{n} 
	\quad \textrm{ with } \quad 
	\sum_{n=1}^{\infty} | u_{n} | < \infty,  \quad n \ge 1 ,
\]
the stability region of the convolution quadrature (\ref{eq:GeneralConvolutionQuadrature}) is
\[
	S_{\alpha} = \Cset \backslash \left\{ \frac{1}{\omega_{\alpha}(\xi)} \, : \: | \xi | \le 1 \right\} 
	, \quad
	\omega_{\alpha}(\xi) = \sum_{n=0}^{\infty} \omega_{n} \xi^{n} .
\]
\end{thm}

We refer to \cite{Garrappa2010_IJCM} and \cite{Lubich1986} to verify that the PI rule and the FLMMs investigated in this paper fulfill the assumptions of Theorem \ref{thm:FED_StabilityRegion}.

In the numerical treatment of ODEs $A$--stability means that the whole complex semi--plane is included in the stability region of the method. $A$--stable methods are therefore particularly attractive since they ensure that the numerical solution of stable linear problems converges to $0$ regardless of the size of the step $h$. Because of the result in Theorem \ref{thm:FDEStability}, for FDEs it is more useful to use the concept of $A(\alpha\frac{\pi}{2})$--stability, which means that the whole sector $\Sigma_{\alpha}$ must be included in the stability region.

%$A$--stable methods are particularly attractive in the numerical treatment of ODEs since they ensure that the numerical solution of stable linear problems converges to $0$ regardless of the size of the step $h$. Because of the result in Theorem \ref{thm:FDEStability}, with FDEs the concept of $A$--stability must be clearly replaced by $A(\alpha\frac{\pi}{2})$--stability, which means that the whole sector $\Sigam_{\alpha}$ must be included in the region of stability of the method.

While $\omega_{\alpha}(\xi)$ is a priori known for FLMMs, the generating function of PI rules does not possess an analytical representation and, therefore, it can be just numerically approximated. 

\begin{figure}[htb]
	\centering
	\begin{tabular}{c@{\hspace{0.3cm}}c@{\hspace{0.3cm}}c@{\hspace{0.3cm}}c}
		PI U & FT  & NG  & FBDF \\
		\includegraphics[width=0.22\textwidth,height=0.23\textwidth]{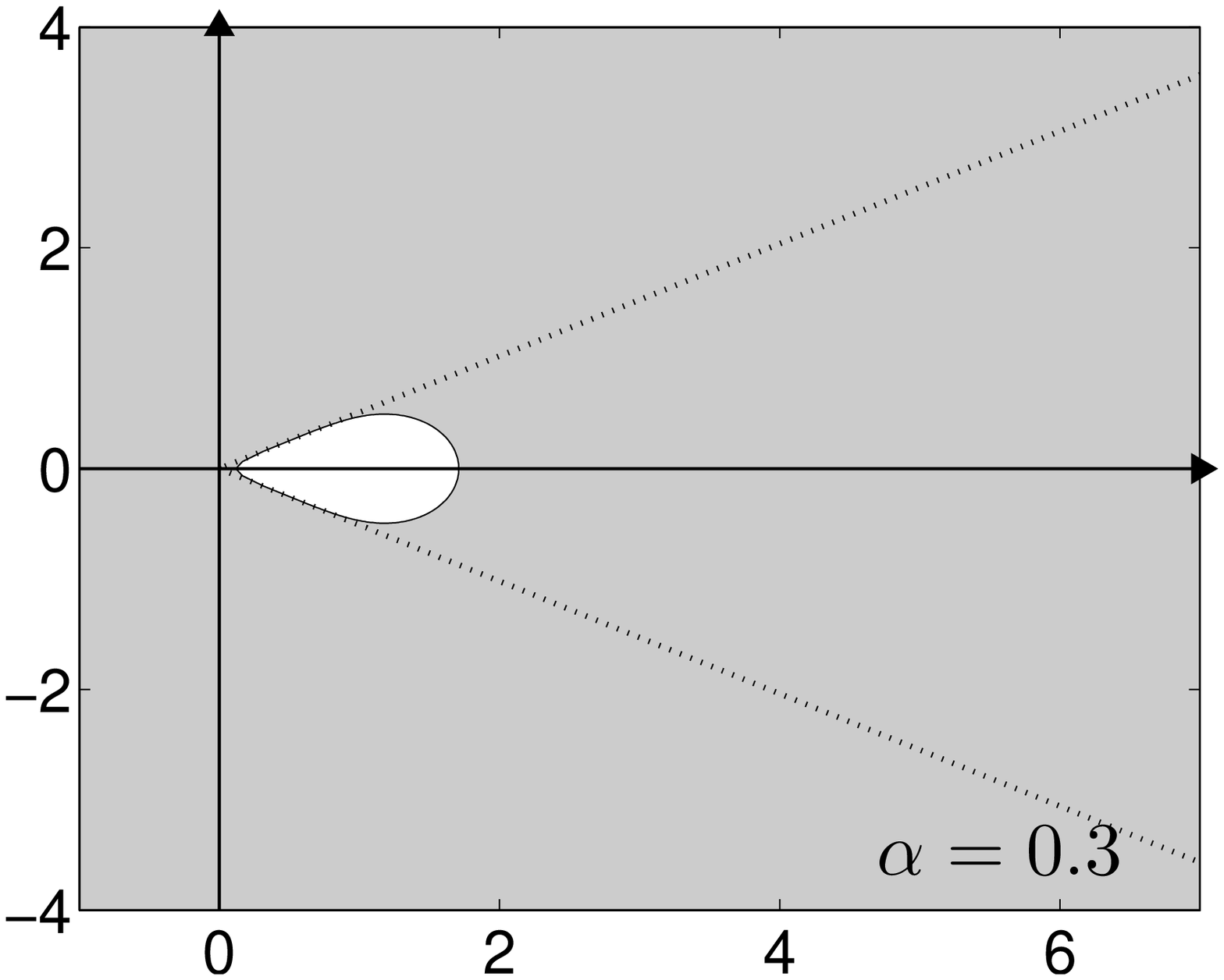} 
		&
		\includegraphics[width=0.22\textwidth,height=0.23\textwidth]{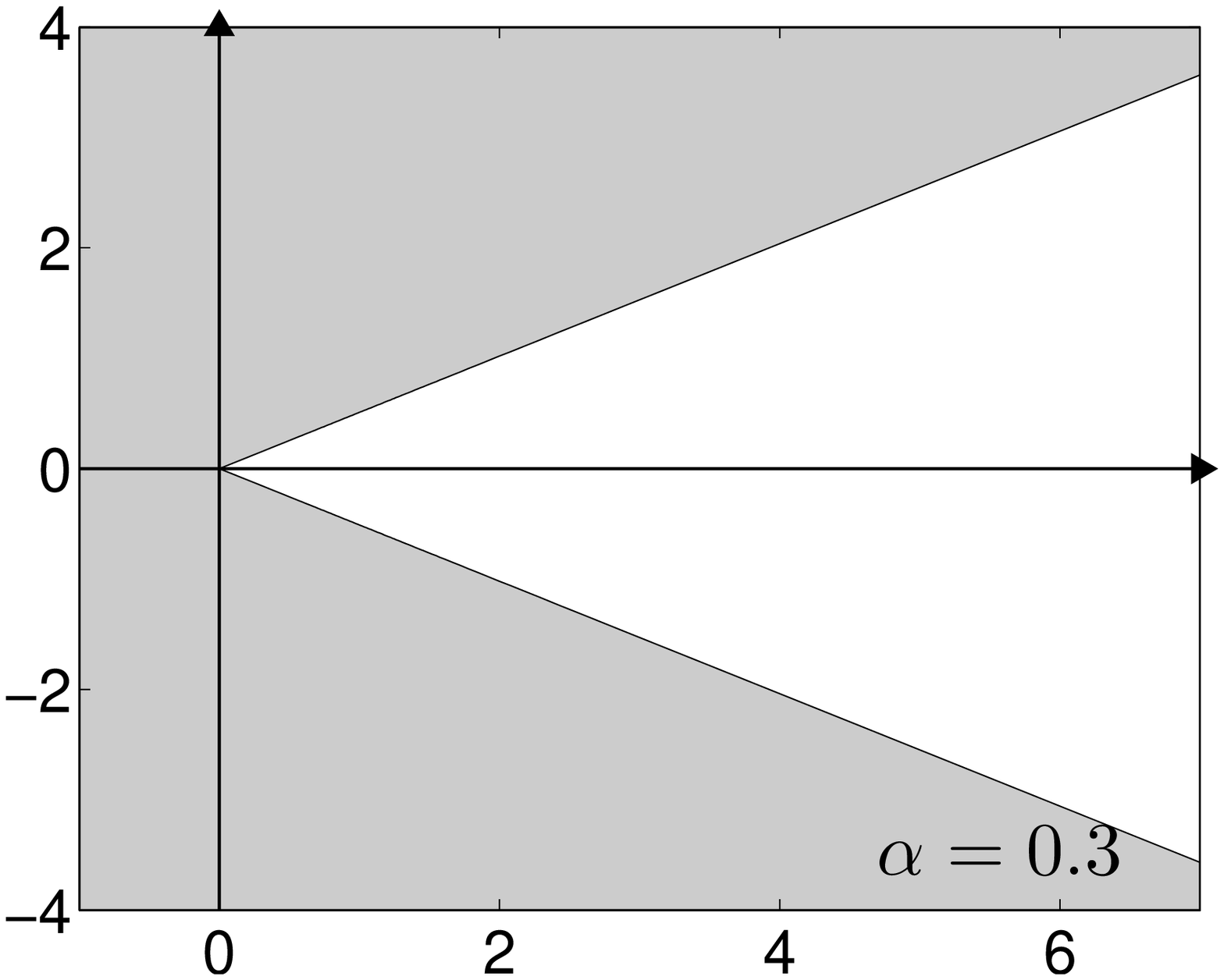} 
		&	
		\includegraphics[width=0.22\textwidth,height=0.23\textwidth]{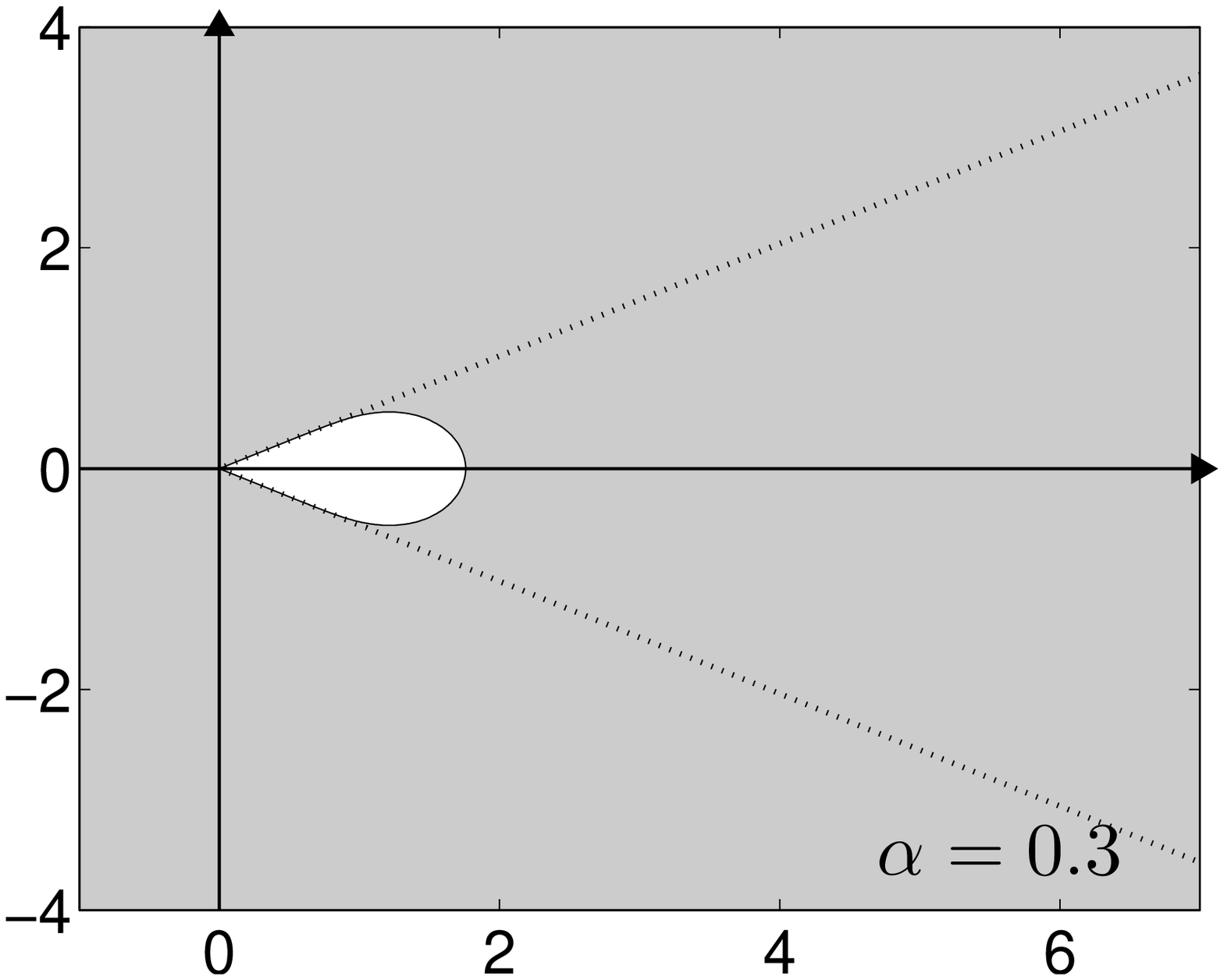} 
		&
		\includegraphics[width=0.22\textwidth,height=0.23\textwidth]{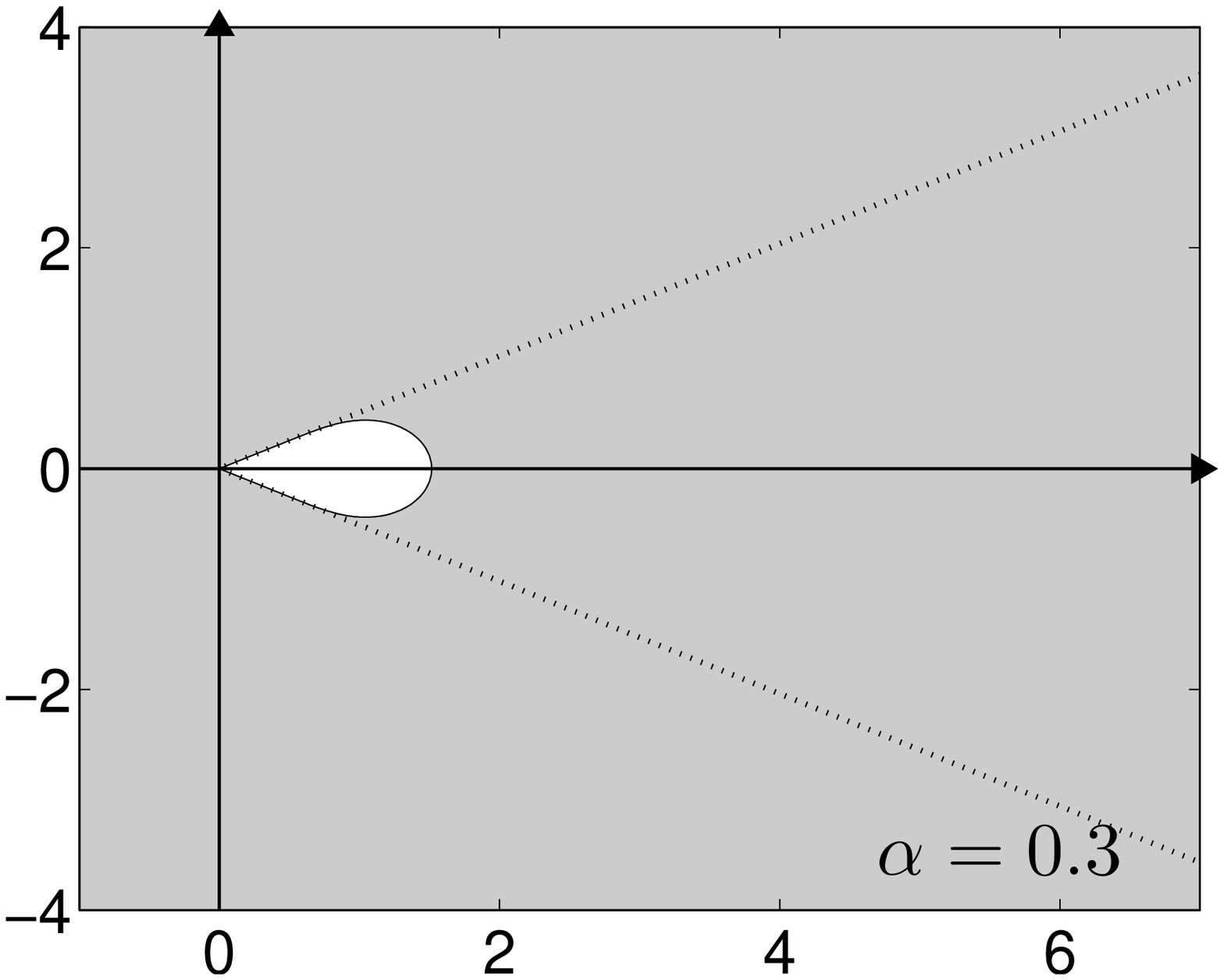} 
		\\	
		\includegraphics[width=0.22\textwidth,height=0.23\textwidth]{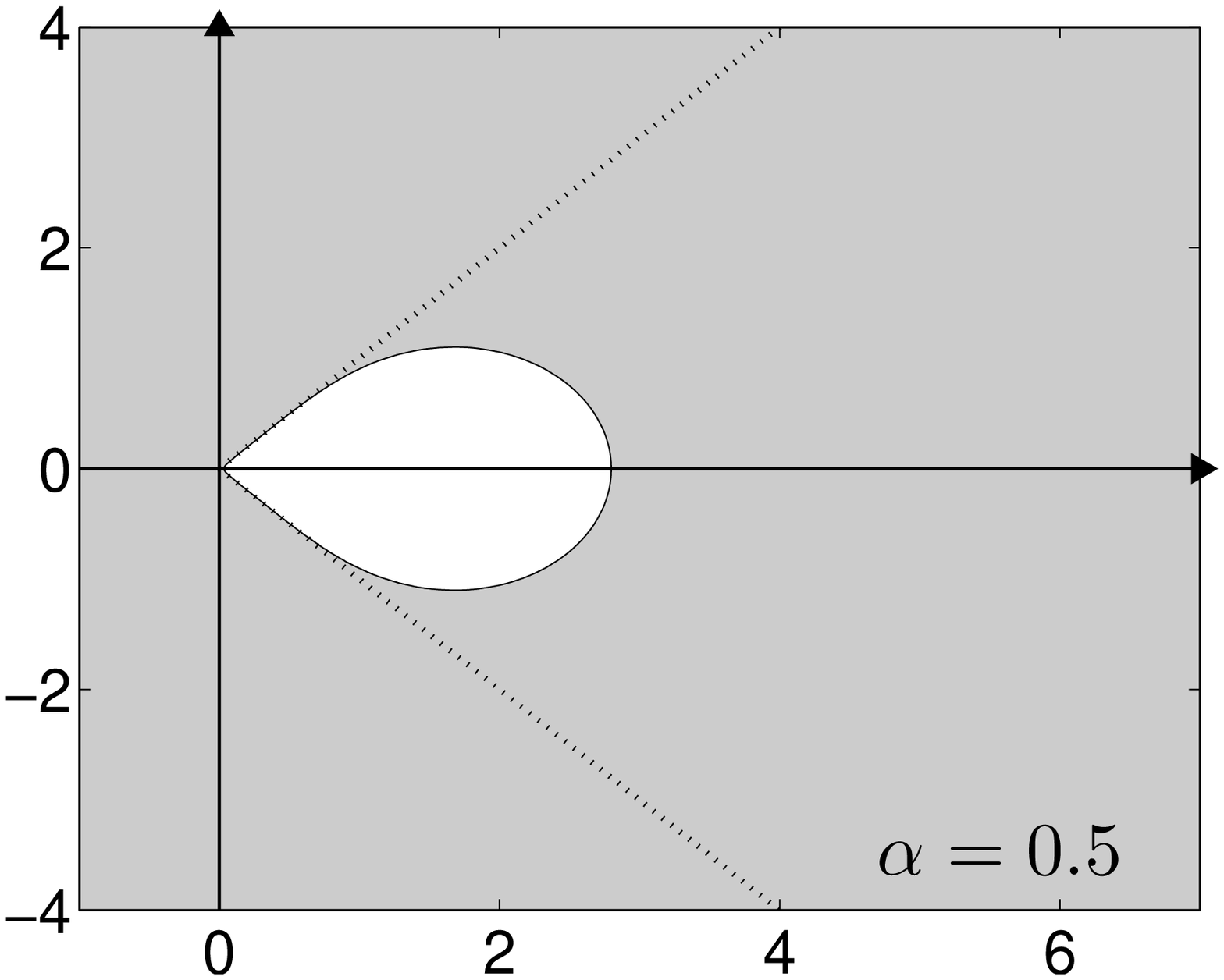} 
		&
		\includegraphics[width=0.22\textwidth,height=0.23\textwidth]{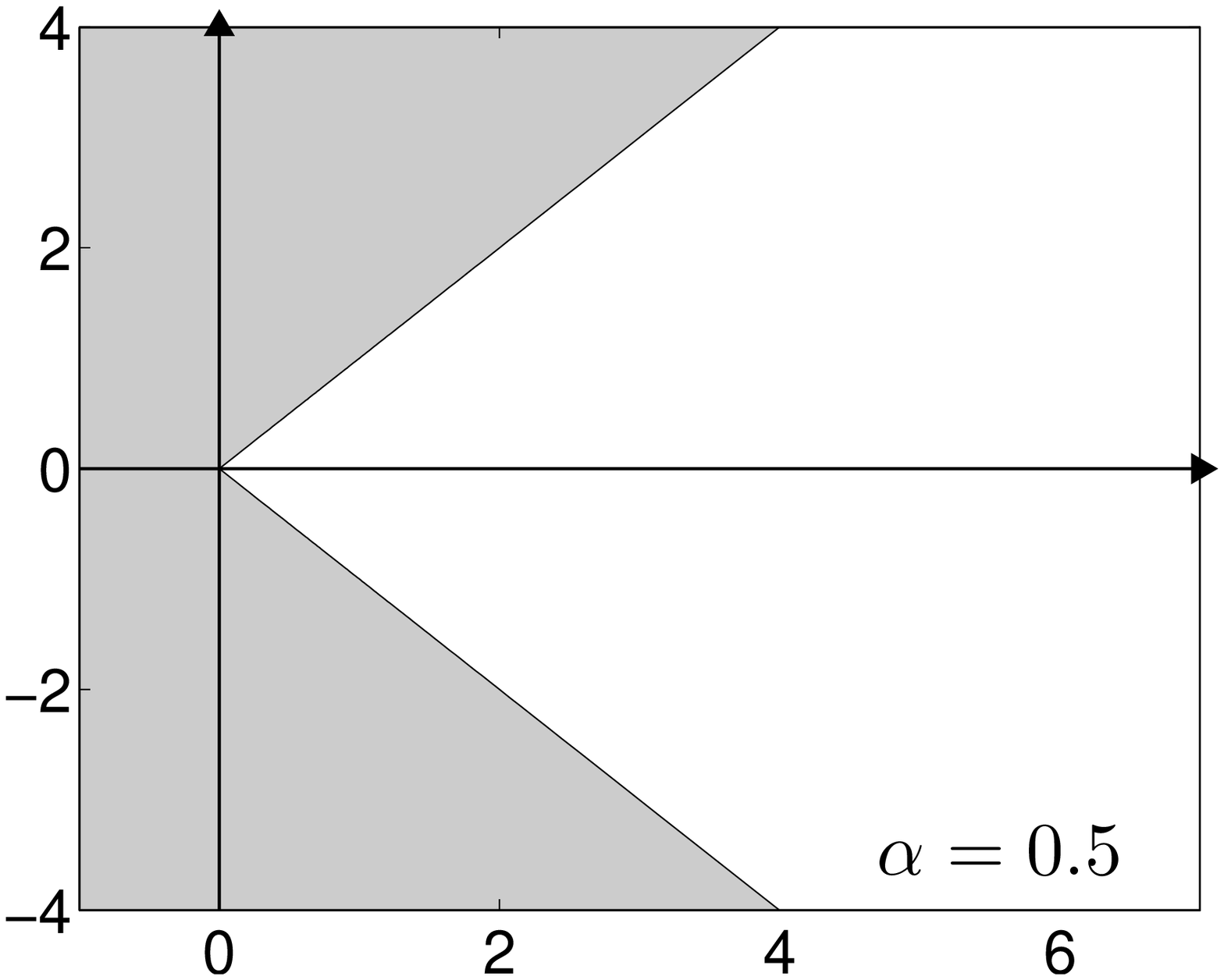} 
		&	
		\includegraphics[width=0.22\textwidth,height=0.23\textwidth]{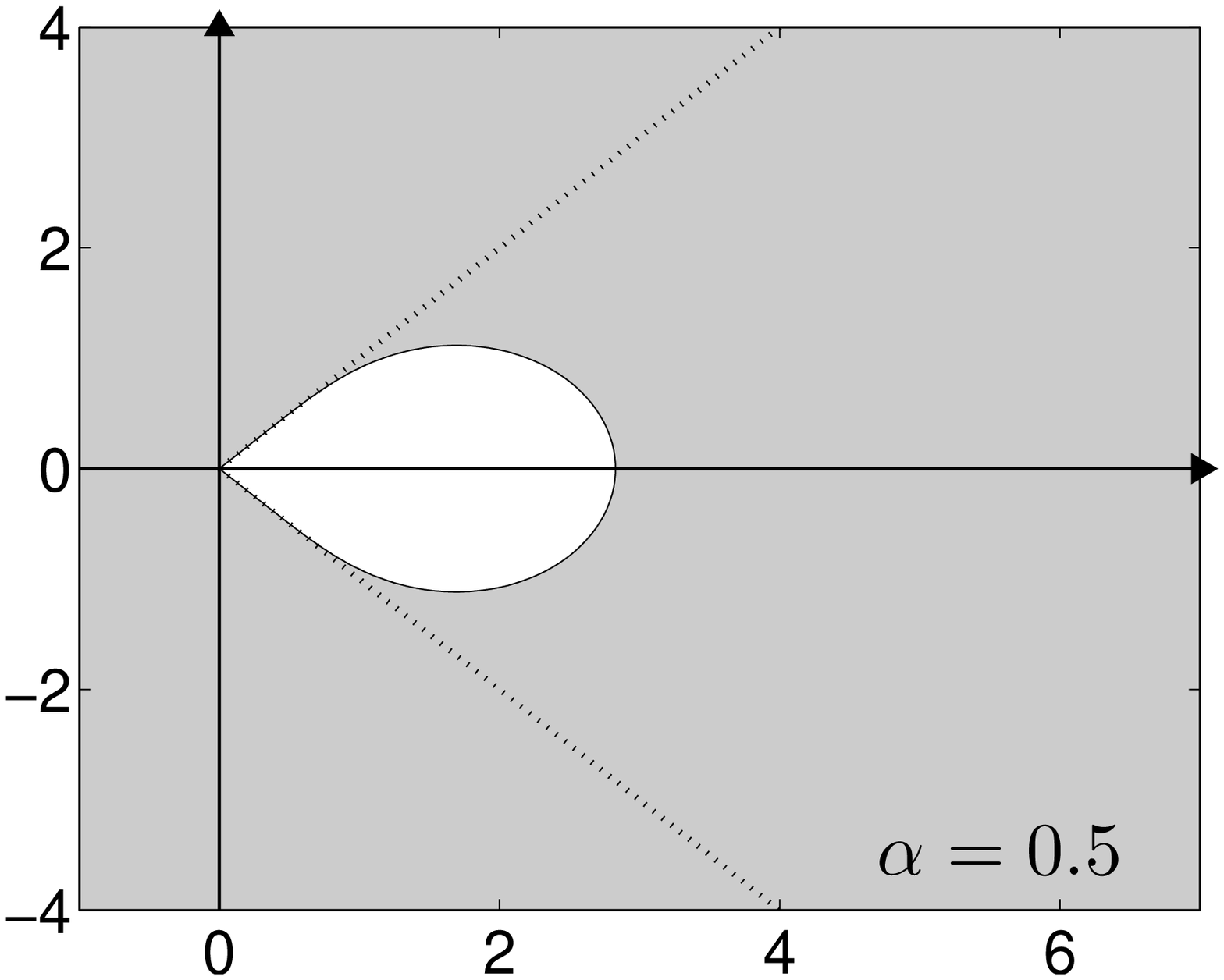} 
		&
		\includegraphics[width=0.22\textwidth,height=0.23\textwidth]{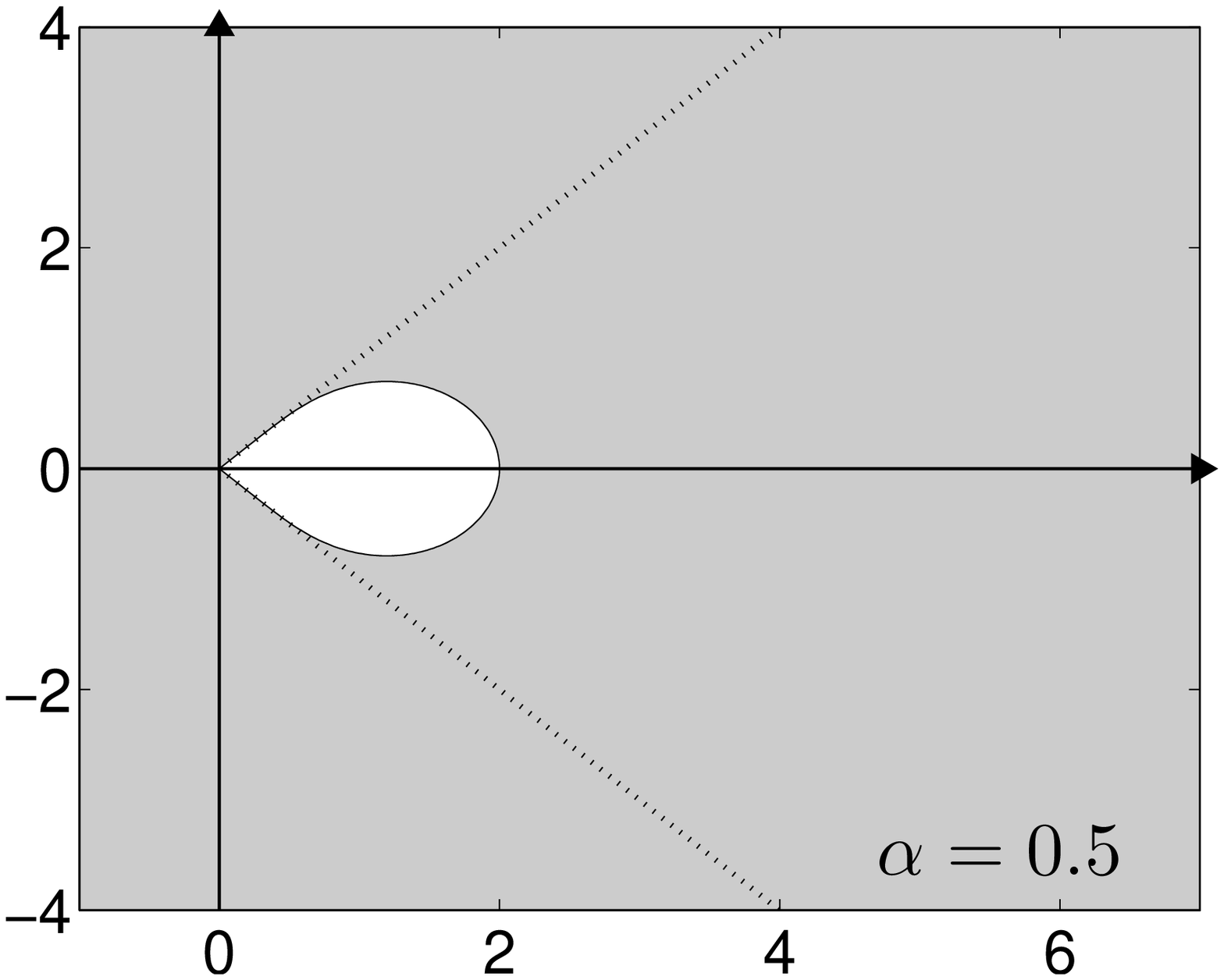} 
		\\
		\includegraphics[width=0.22\textwidth,height=0.23\textwidth]{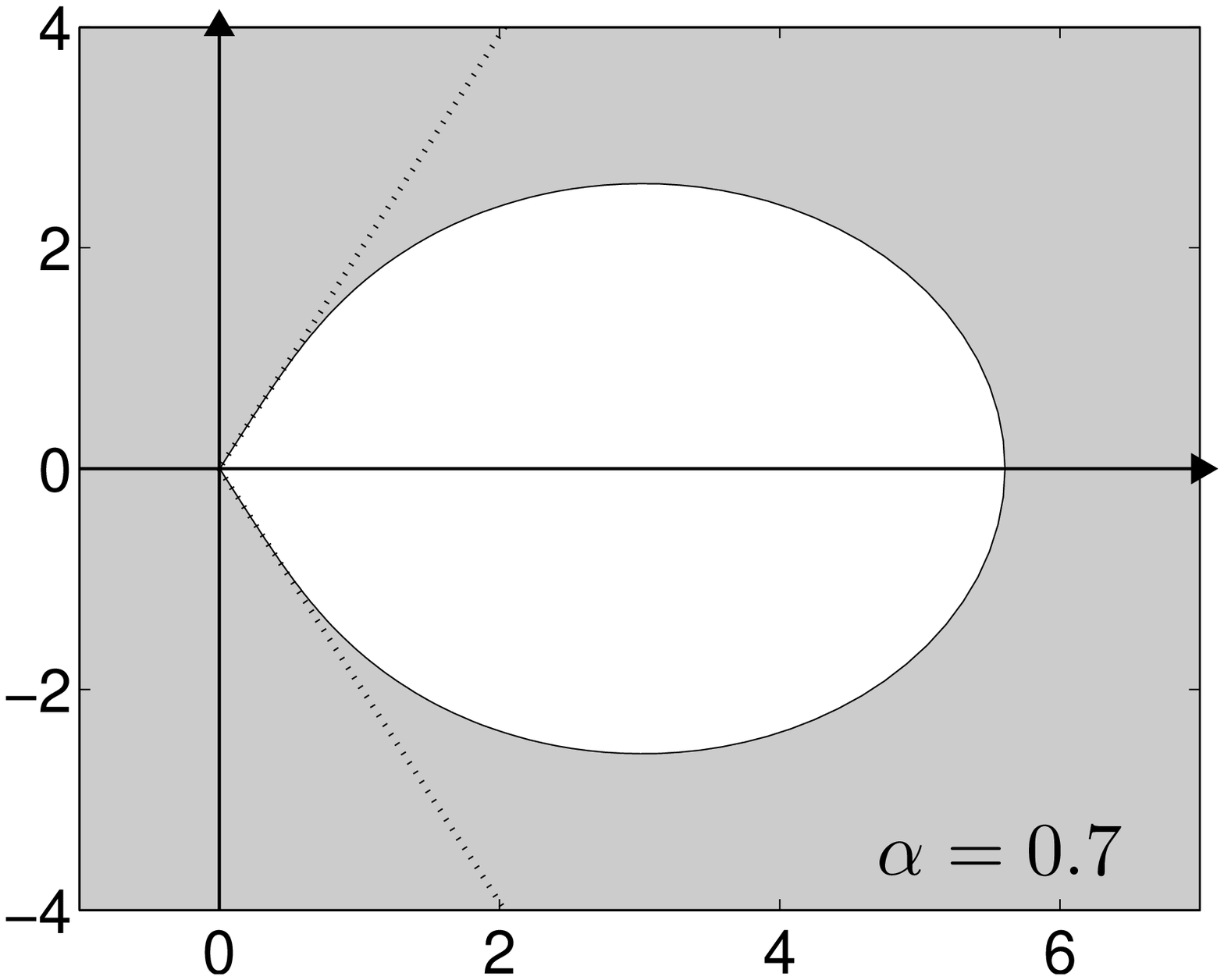}
		&
		\includegraphics[width=0.22\textwidth,height=0.23\textwidth]{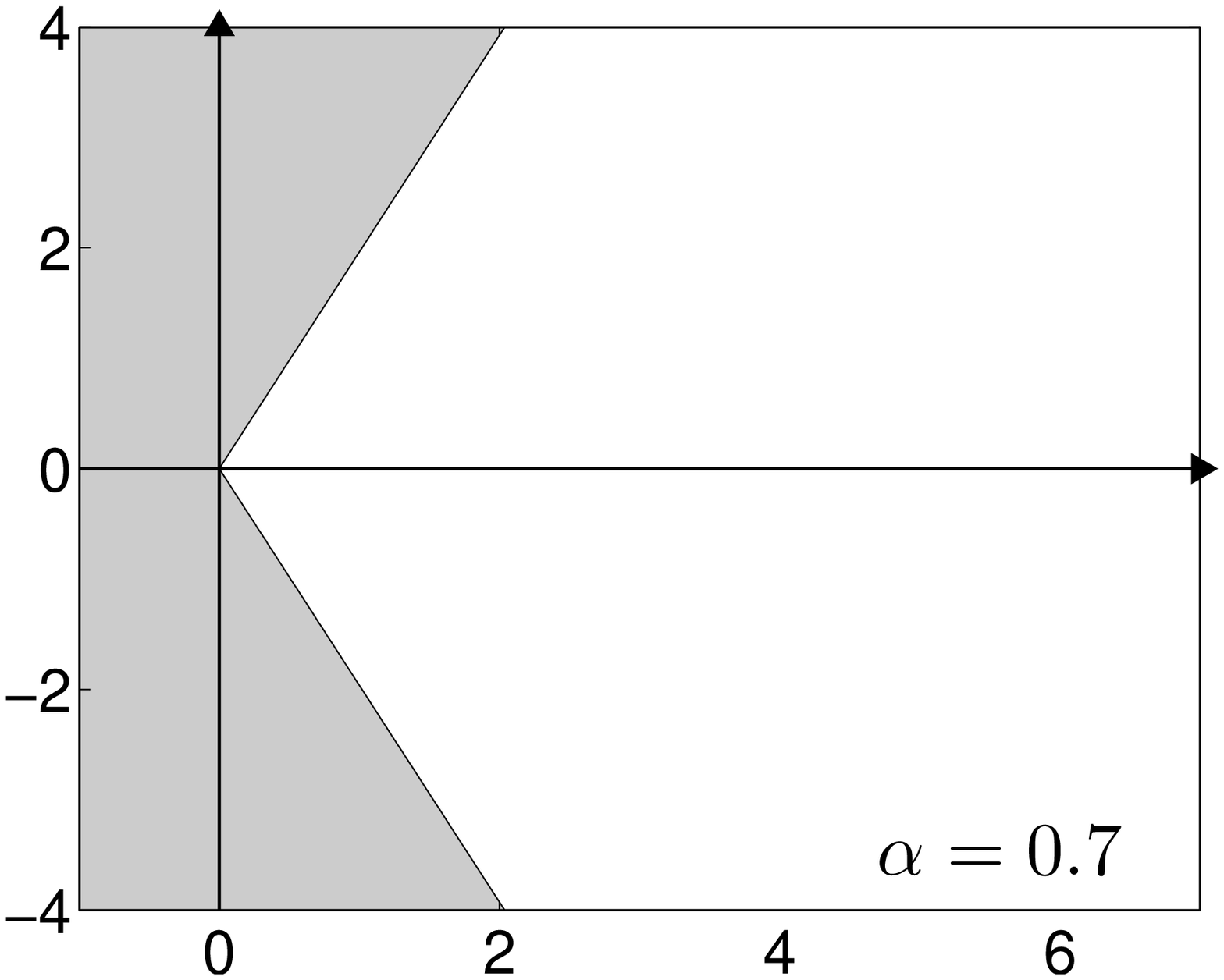} 
		&	
		\includegraphics[width=0.22\textwidth,height=0.23\textwidth]{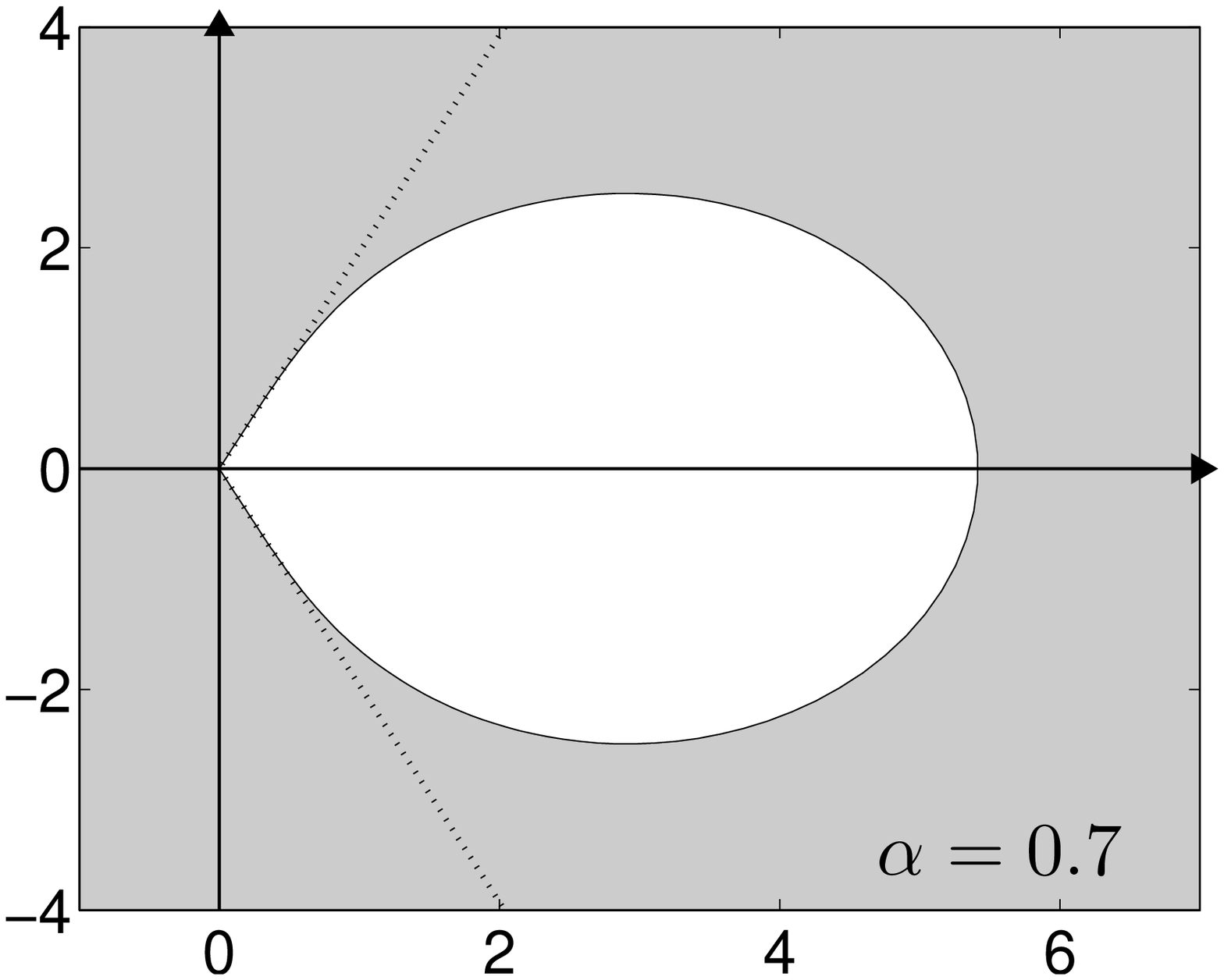} 
		&
		\includegraphics[width=0.22\textwidth,height=0.23\textwidth]{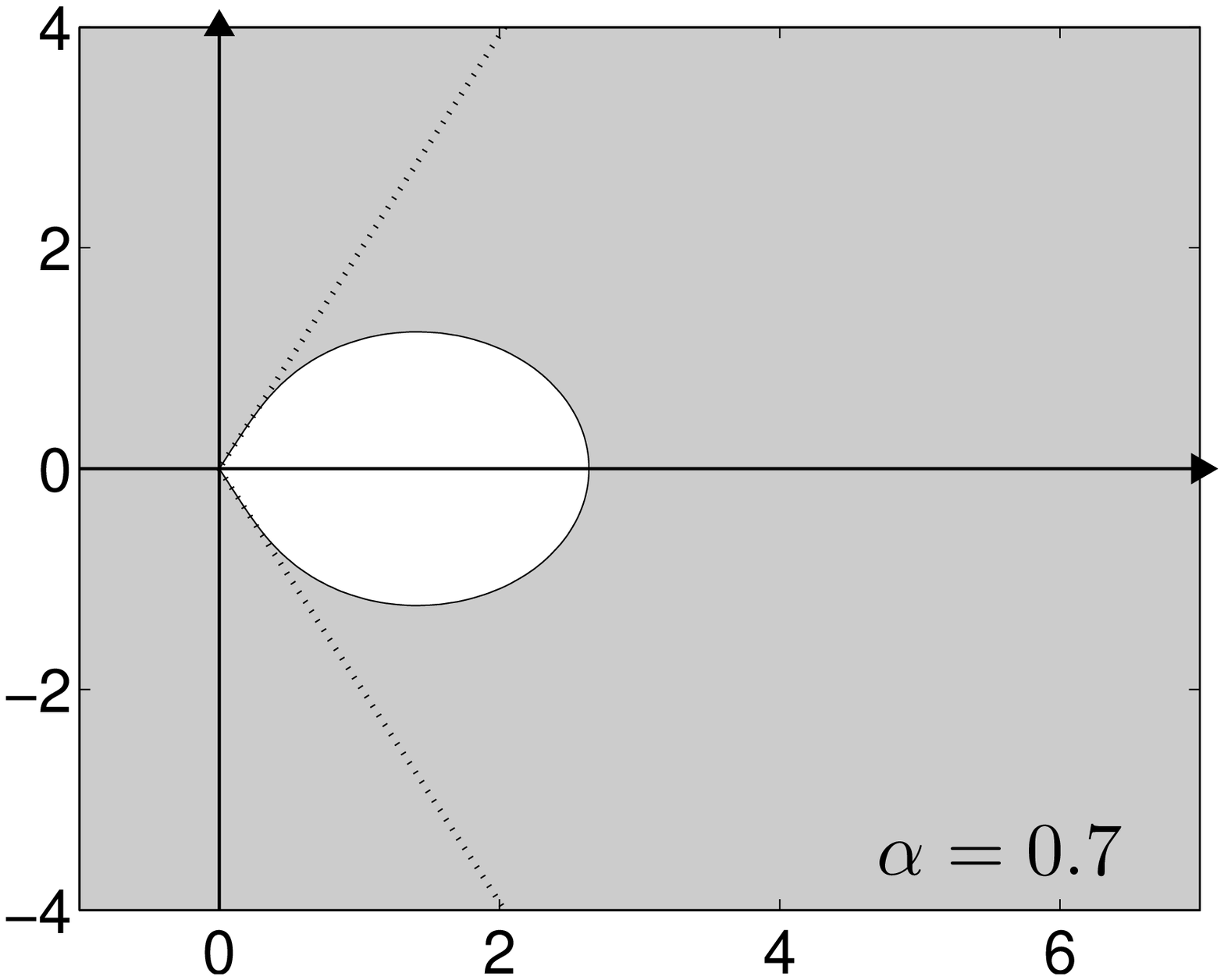} 
		\\
	\end{tabular}
\caption{Stability regions (gray areas) for some orders $0<\alpha<1$}
\label{fig:Fig_StabilityRegions1}
\end{figure}

In Figure \ref{fig:Fig_StabilityRegions1} we present a comparison of the stability regions (the gray areas) of the methods discussed in the previous sections for various values $0<\alpha<1$ (dotted lines denote the corresponding sectors of stability $\Sigma_{\alpha}$). To identify each method, here and throughout the following section we make use of the acronyms listed in Table \ref{T:Acronyms}.

\begin{table}
\centering
\begin{tabular}{l|l|c} \hline
 Acronym & Formula & Subsection \\ \hline
	PI U & PI rule on uniform grid & \ref{SS:PI_Uniform} \\
	PI G & PI rule on graded grid & \ref{SS:PI_Graded} \\ 
	FT   & fractional trapezoidal rule & \ref{SS:FracTrapRule} \\
	NG   & Newton--Gregory formula & \ref{SS:NewtonGregoryFormula} \\
	FBDF &  fractional BDF formula & \ref{SS:FracBDF} \\ \hline	 
\end{tabular}
\caption{Acronyms used to denote methods}\label{T:Acronyms}
\end{table}

We can distinctly observe that, as $\alpha$ increases from $0$ to $1$, the stability region of each method decreases but always includes $\Sigma_{\alpha}$; we can draw the conclusion that all the methods are $A(\alpha\frac{\pi}{2})$--stable for $0<\alpha<1$.

Within ODEs, the stability region of the trapezoidal rule is the negative semi--plane $\Cset^{-}$ and exactly corresponds to the region in which the true solution converges to the steady--state. Among the methods presented in this paper, just FT preserves this feature and indeed its stability regions coincide with $\Sigma_{\alpha}$ for any value of $\alpha$. NG and PI have regions of almost the same size but larger with respect to FT; as in the ODE case, FBDF possesses the largest regions of stability.

%It is well known that the stability region of the trapezoidal rule for ODEs is the negative semi--plane $\Cset^{-}$, which exactly corresponds to the region in which the true solution converges to the steady--state. Among the methods presented in this paper, just the trapezoidal FLMM preserves this features and indeed its stability regions coincide with $\Sigma_{\alpha}$ for any value of $\alpha$. The Newton Gregory FLMM  and the PI rule have regions of almost the same size but larger with respect to the trapezoidal FLMM. The regions for the FBDF2 are the largest one among the various methods: this is a feature inherited from the ODE case, for which BDFs are the methods with the largest regions of stability. 

\begin{figure}[htb]
	\centering
	\begin{tabular}{c@{\hspace{0.3cm}}c@{\hspace{0.3cm}}c@{\hspace{0.3cm}}c}
		PI U & FT  & NG  & FBDF \\
		\includegraphics[width=0.22\textwidth,height=0.23\textwidth]{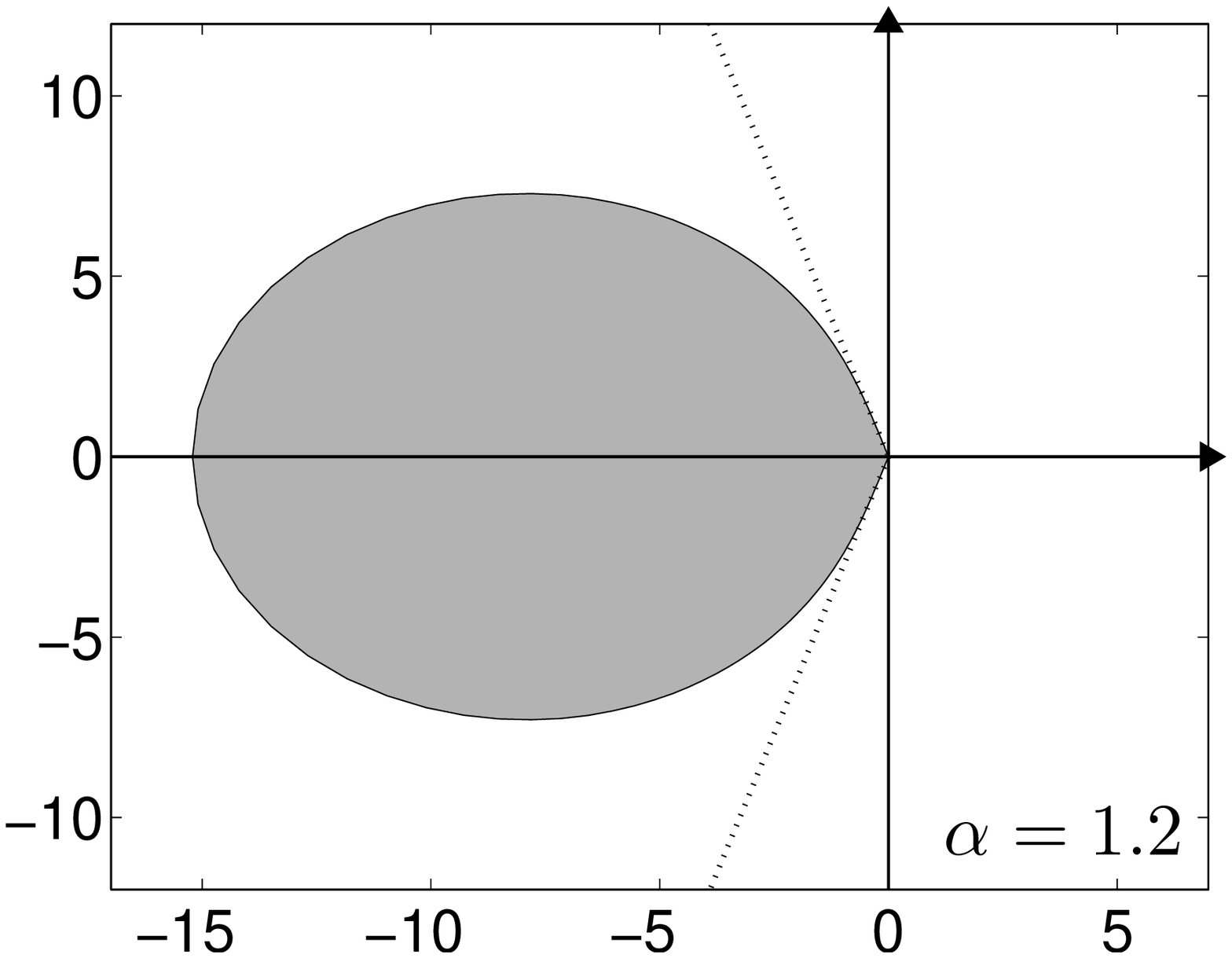} 
		&
		\includegraphics[width=0.22\textwidth,height=0.23\textwidth]{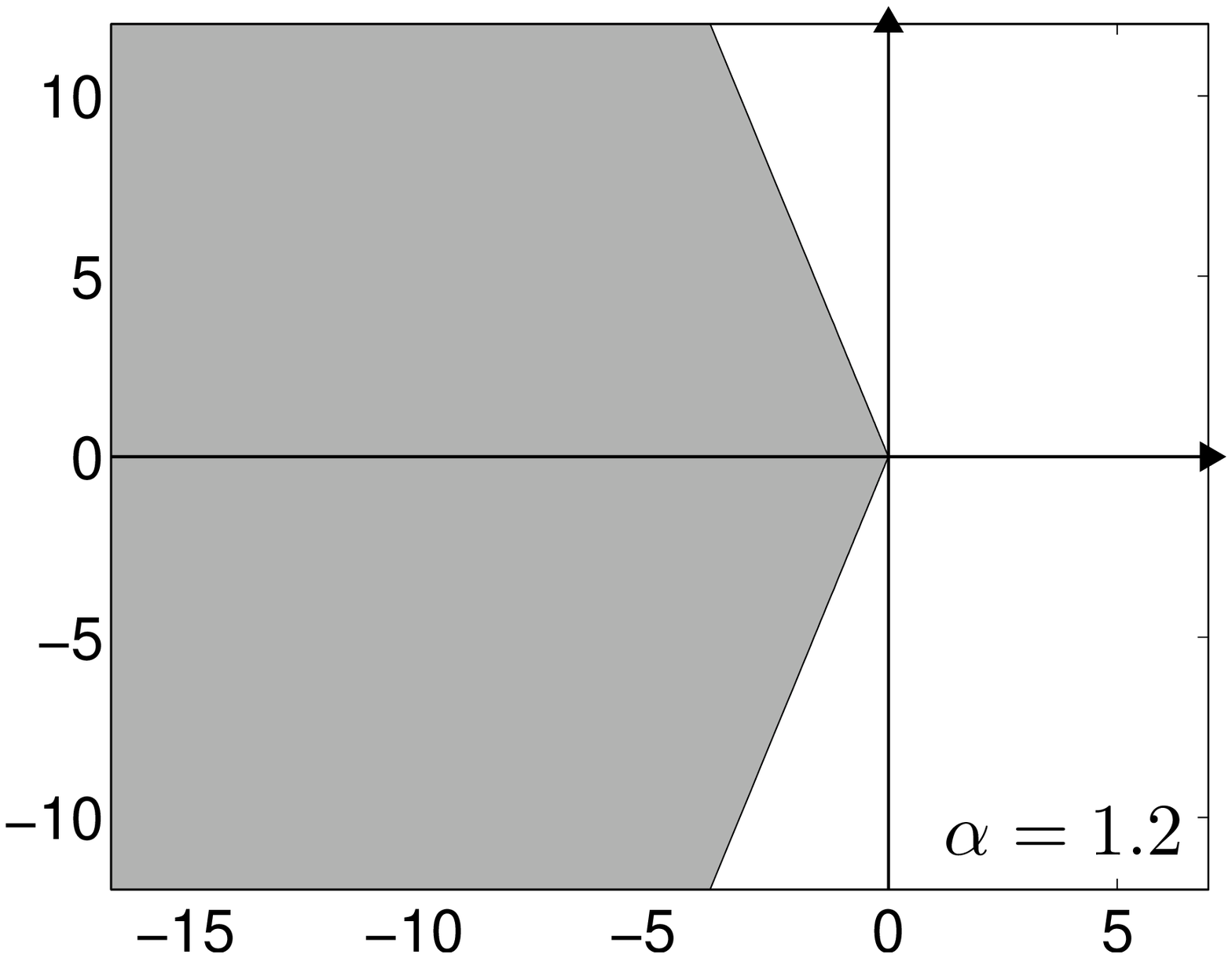} 
		&	
		\includegraphics[width=0.22\textwidth,height=0.23\textwidth]{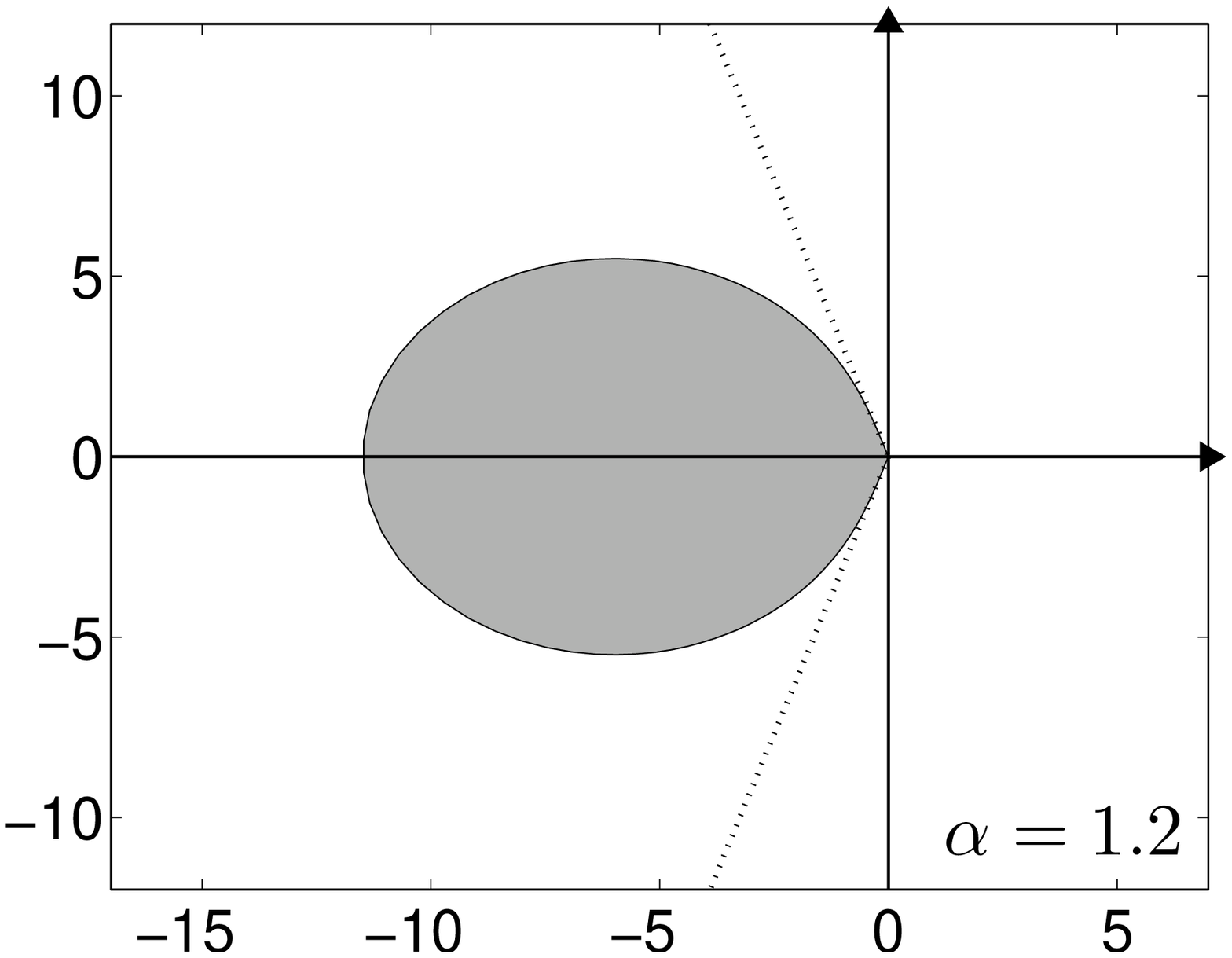} 
		&
		\includegraphics[width=0.22\textwidth,height=0.23\textwidth]{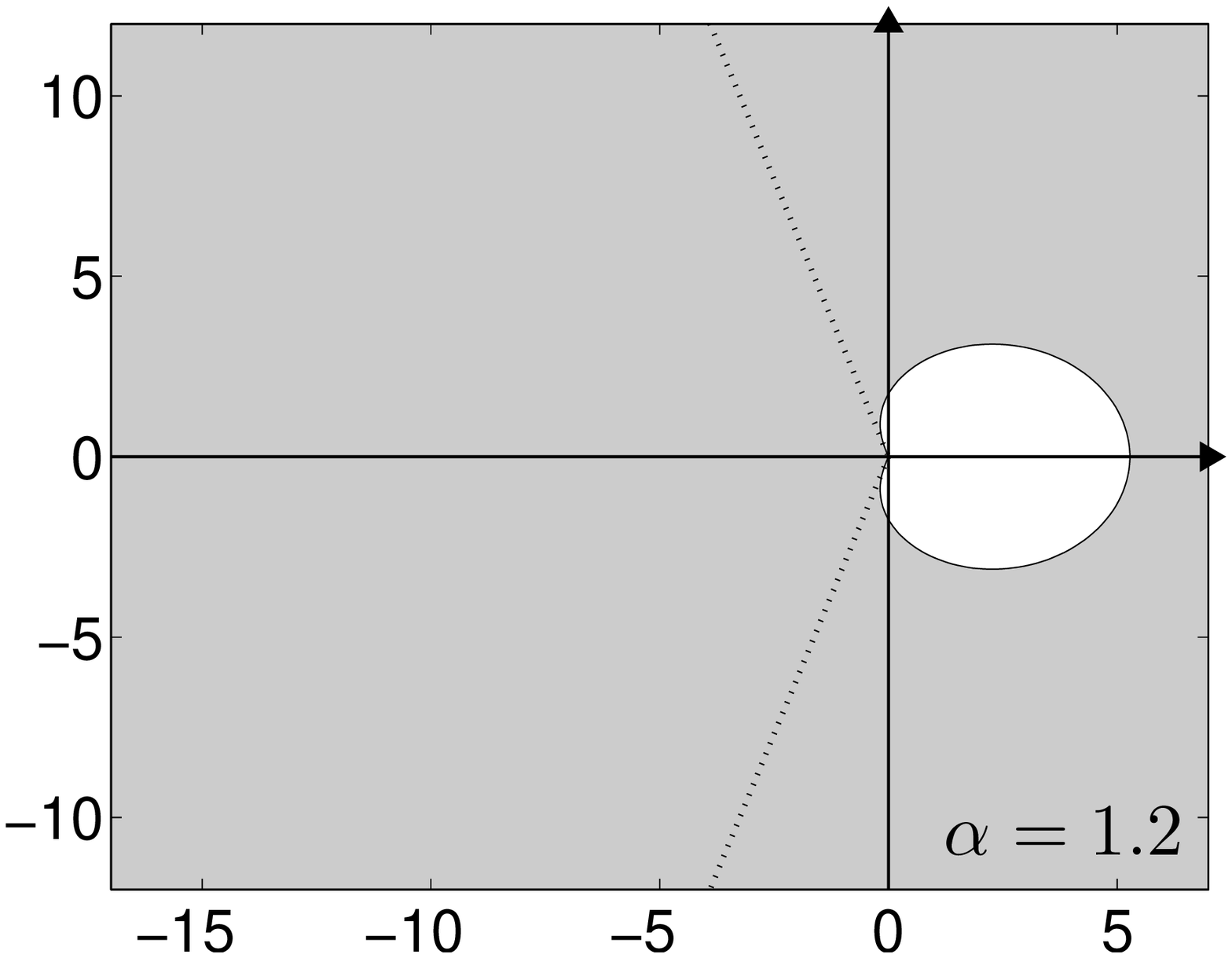} 
		\\	
		\includegraphics[width=0.22\textwidth,height=0.23\textwidth]{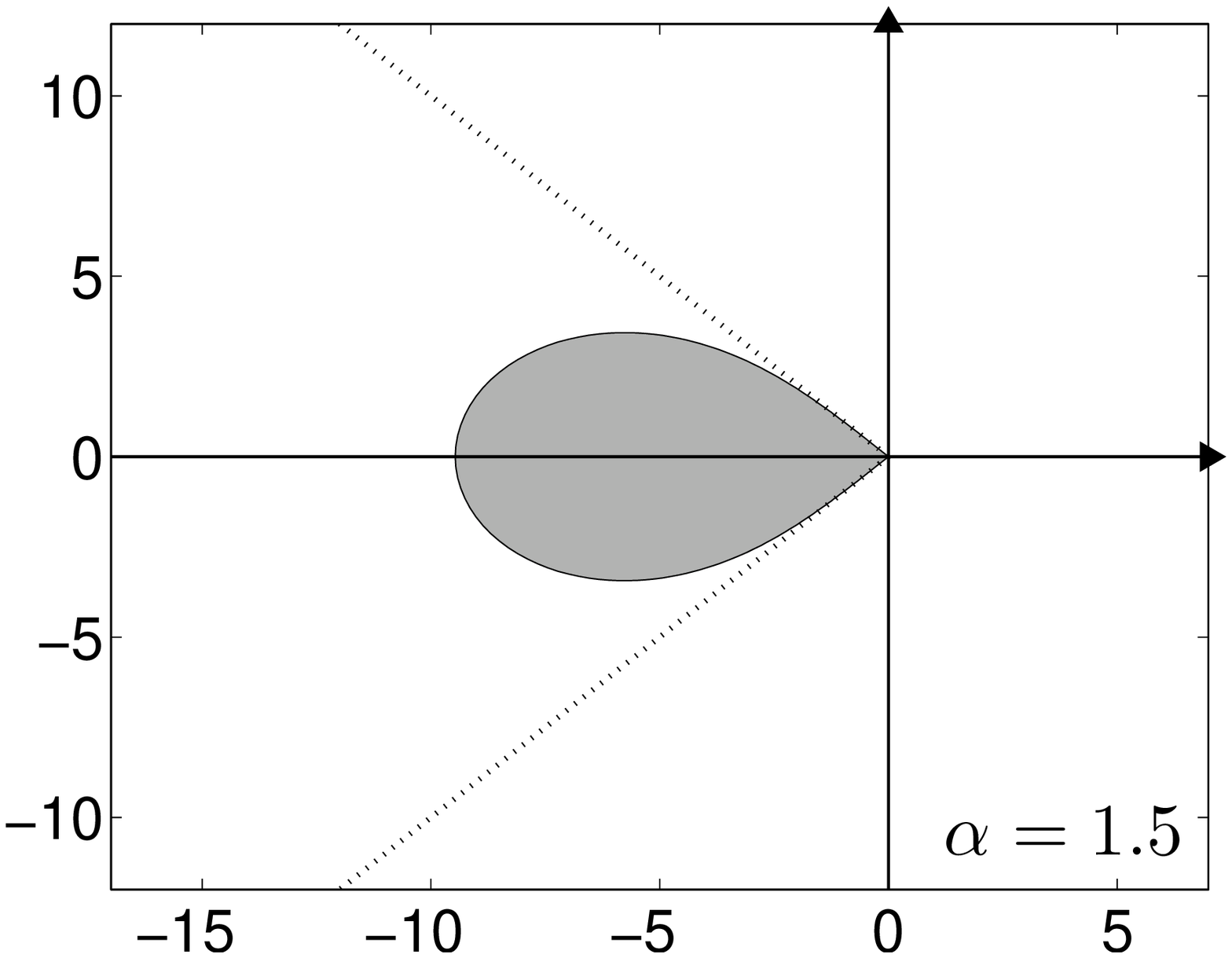} 
		&
		\includegraphics[width=0.22\textwidth,height=0.23\textwidth]{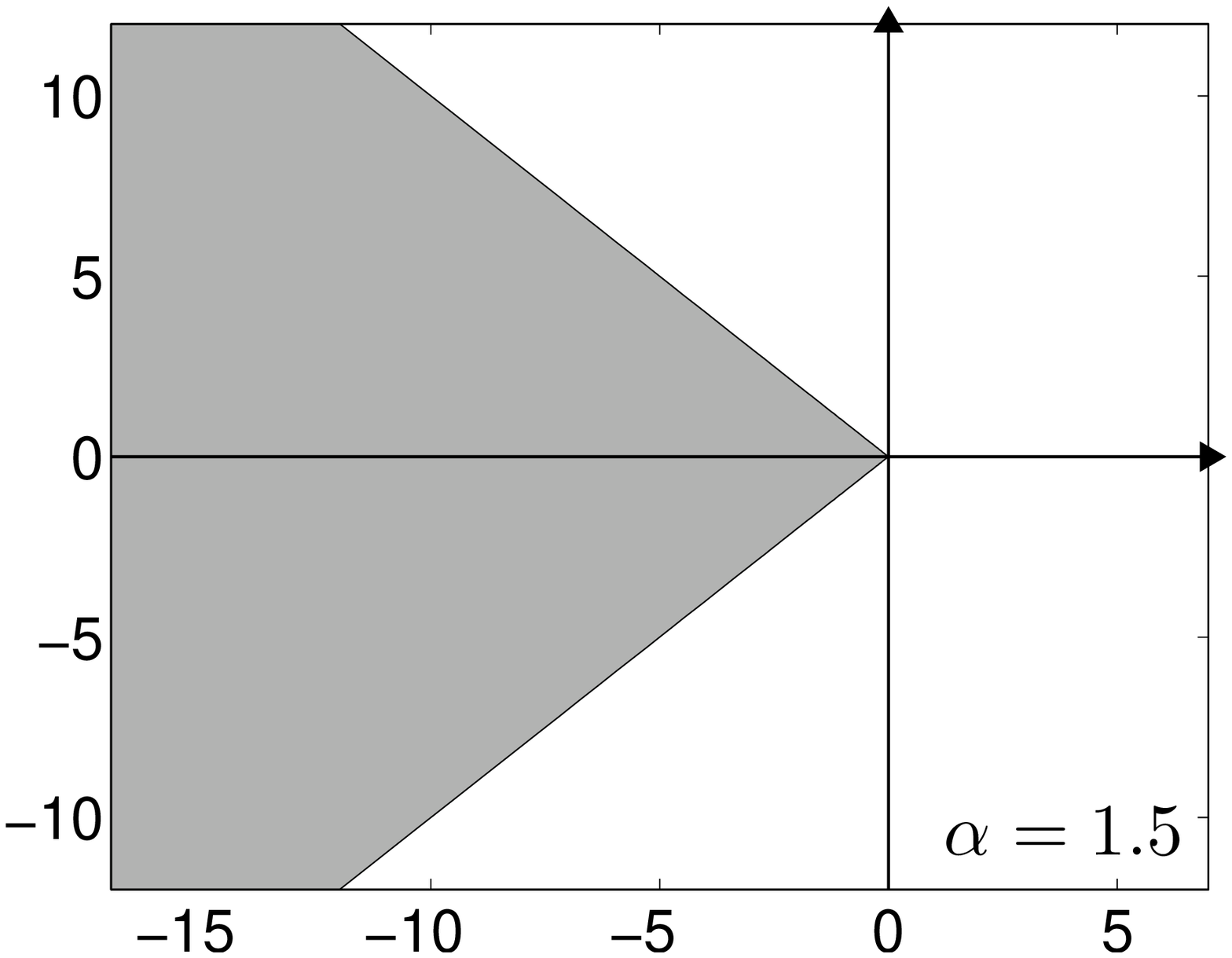} 
		&	
		\includegraphics[width=0.22\textwidth,height=0.23\textwidth]{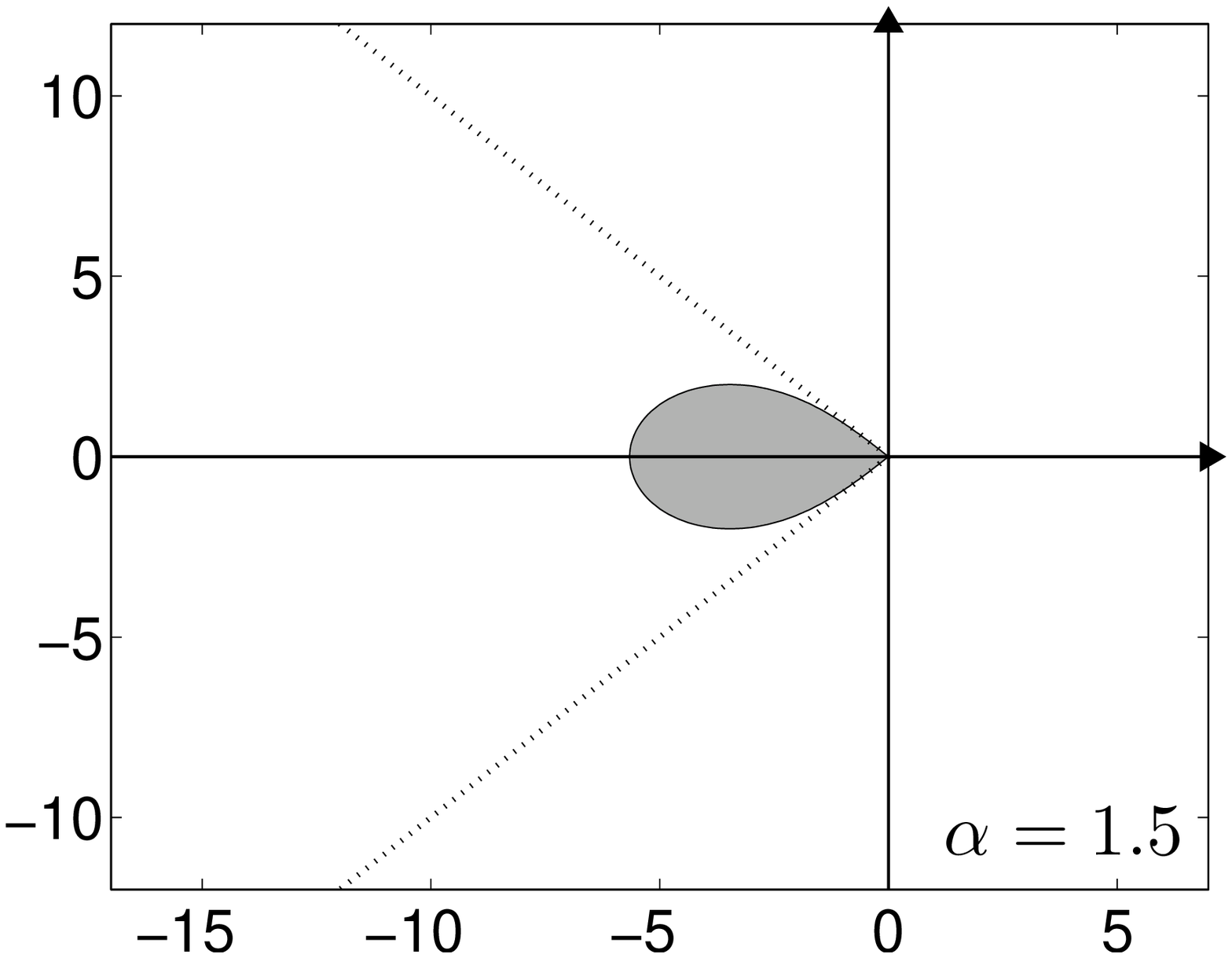} 
		&
		\includegraphics[width=0.22\textwidth,height=0.23\textwidth]{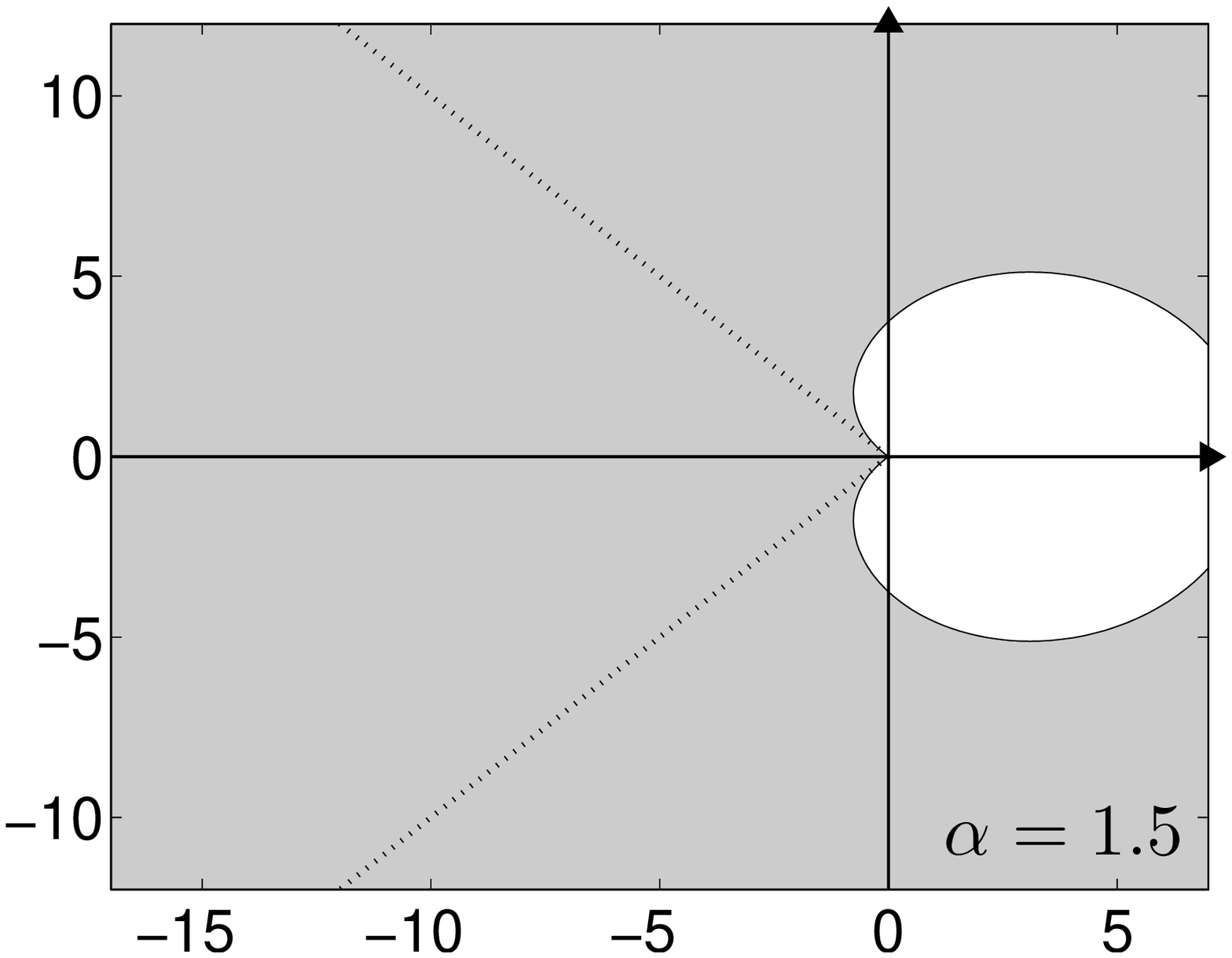} 
		\\
		\includegraphics[width=0.22\textwidth,height=0.23\textwidth]{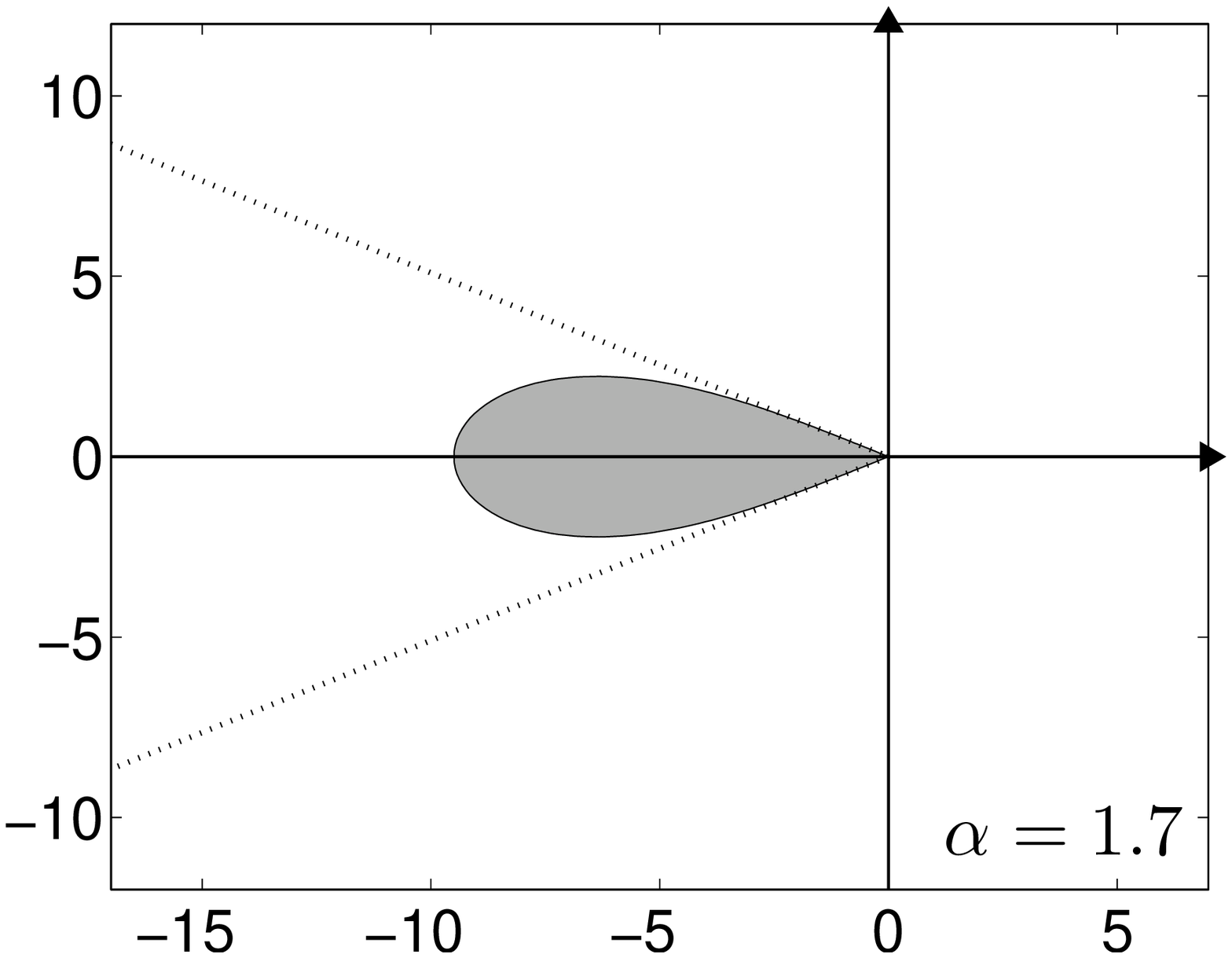} 
		&
		\includegraphics[width=0.22\textwidth,height=0.23\textwidth]{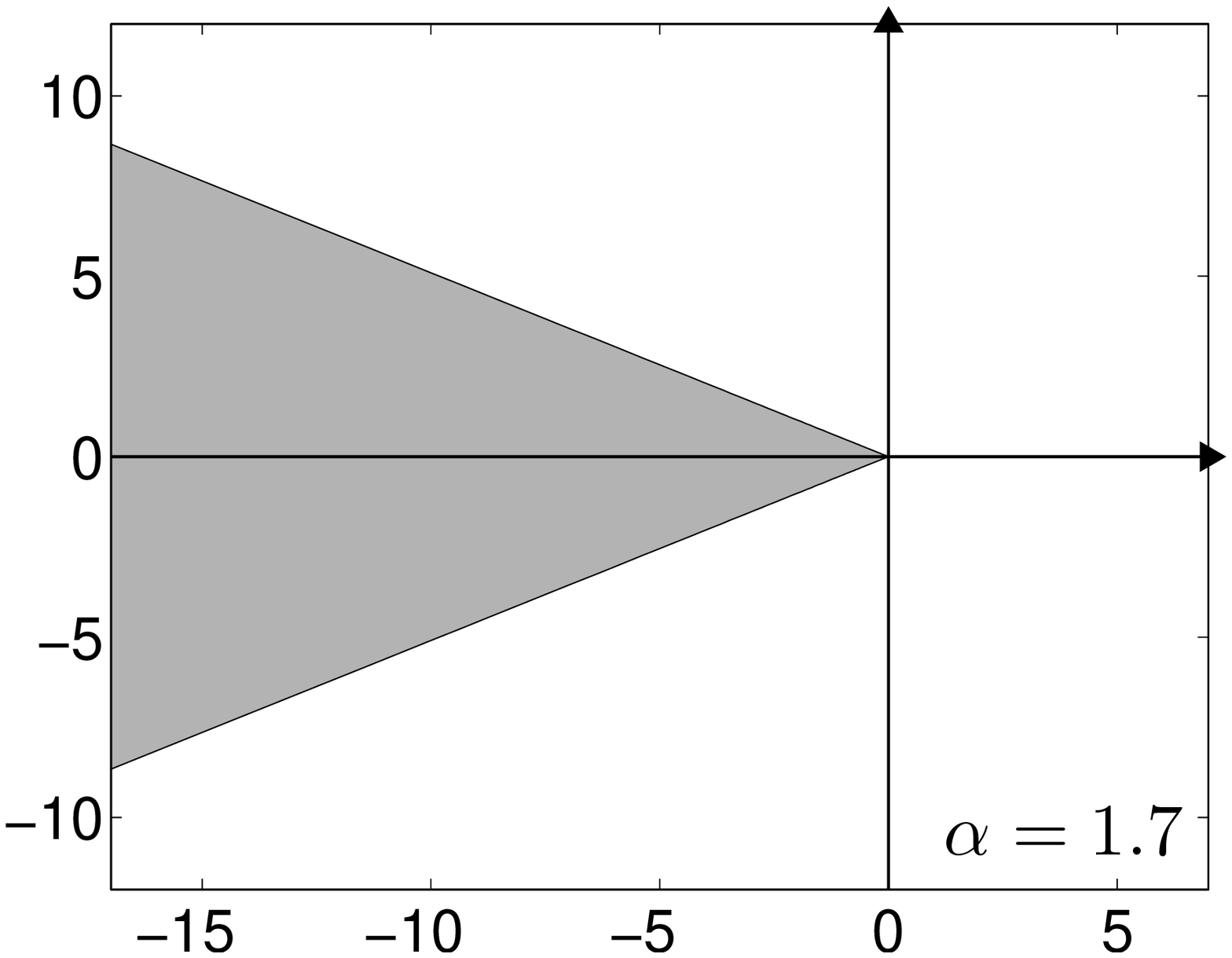} 
		&	
		\includegraphics[width=0.22\textwidth,height=0.23\textwidth]{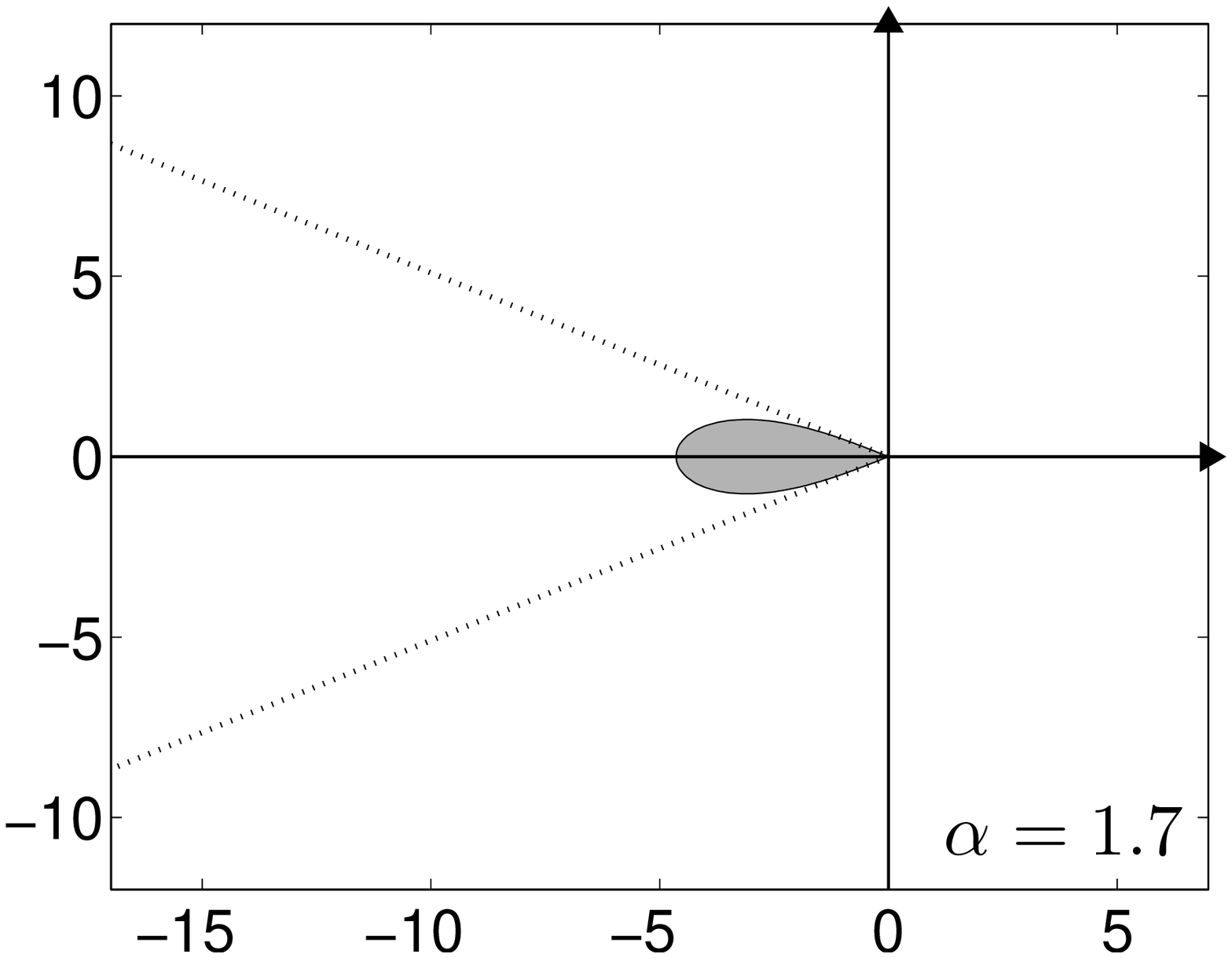} 
		&
		\includegraphics[width=0.22\textwidth,height=0.23\textwidth]{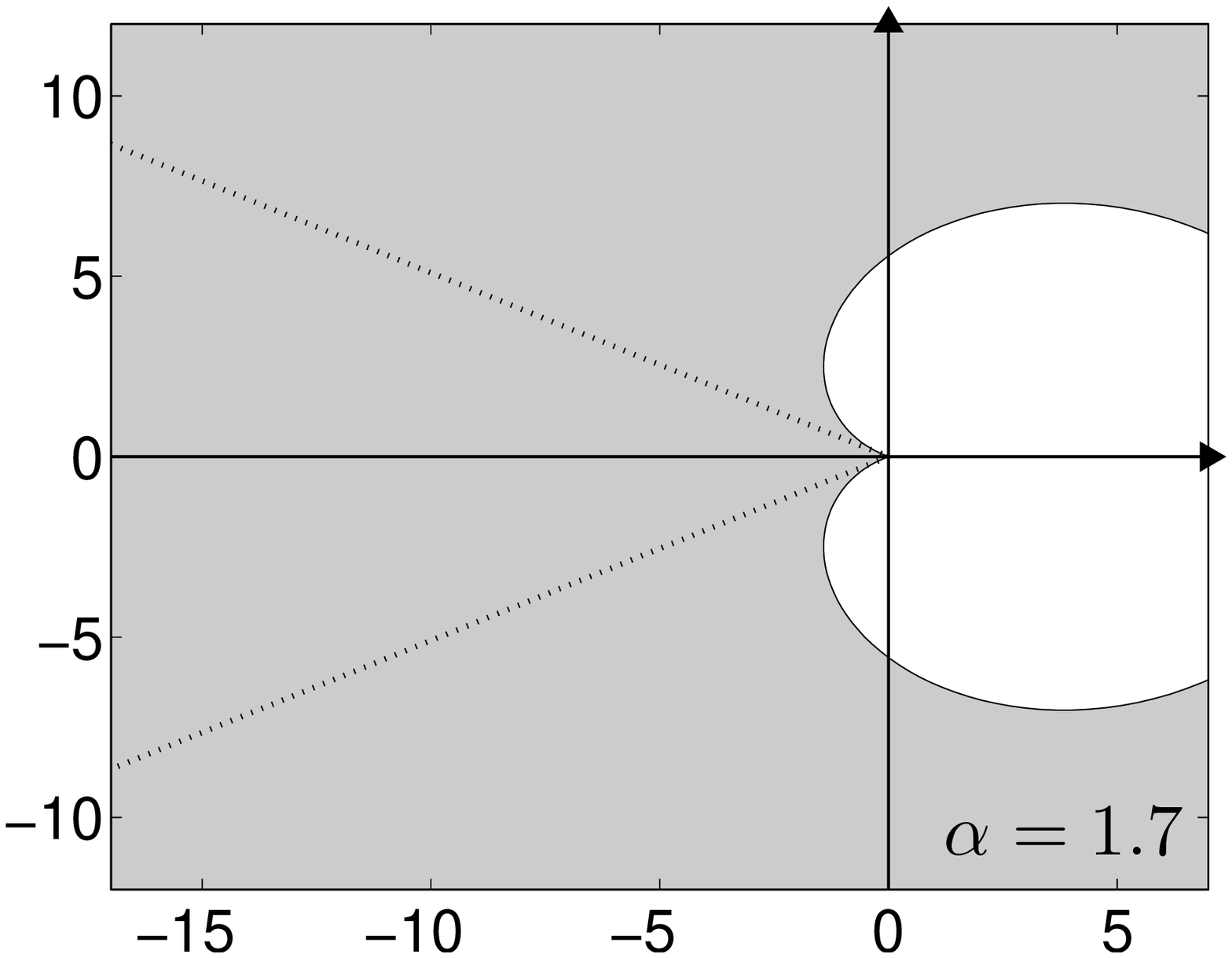} 
		\\
	\end{tabular}
\caption{Stability regions (gray areas) for some orders $\alpha>1$}
\label{fig:Fig_StabilityRegions2}
\end{figure}

Interestingly, the situation is completely different for $1<\alpha <2$ as shown in Figure \ref{fig:Fig_StabilityRegions2}. Whilst the stability regions of FT still coincide with the stability sectors $\Sigma_{\alpha}$ and the stability properties of FBDF still remain the most remarkable ones, we observe that NG and PI completely lose the $A(\alpha\frac{\pi}{2})$--stability. This is an important aspect to be taken into account: despite the implicit nature of these methods, numerical instability can arise when a not sufficiently small step--size is used and, therefore, they seem not suited to solve stiff problems when $1<\alpha<2$.

\section{Implementation details and Matlab codes}\label{S:ImpementationIssues}

In this section we discuss some implementation issues that are essential for the development of efficient codes.

\subsection{Evaluation of the lag term}

The evaluation of the (usually long) lag term is one of the most expensive tasks in convolution quadratures like (\ref{eq:PI_General_Uniform}) or (\ref{eq:FLMM}). The integration along a large number $N$ of steps is usually a critical bottleneck and, if not well designed, it can be an impracticable task involving a number of operations proportional to $N^2$.

Several authors have discussed this problem (see, for instance, \cite{ConteDelPrete2006}); in \cite{Diethelm2011} a strategy based on parallel algorithms has been described. 

In our codes we implement the FFT algorithm described in \cite{HairerLubichSchlichte1985} by which the computational effort downgrades from $N^2$ to $N\log(N)^2$. We do not further discuss technical details of this algorithm, for which we just refer to the above mentioned paper; we simply observe that, without efficient algorithms of this kind, the numerical solution of FDEs very quickly degenerates in an unworkable task as the number $N$ of grid--points increases.

\subsection{Starting weights}

Starting weights $w_{n,j}$ are introduced in (\ref{eq:FLMM}) to cope with the singular behavior of the solution close to the left endpoint of the integration interval. They are evaluated after imposing that (\ref{eq:ConvolutionQuadrature}) is exact for $g(t)=t^{\nu}$ with $\nu \in {\cal A}_{p-1}\cup \bigl\{p-1\bigr\}$ and $p$ the order of convergence of the FLMM \cite{Lubich1986}. 

Because of (\ref{eq:RL_Integral_Power}), and since the methods investigated in this paper have order $p=2$, it is immediate to see that weights $w_{n,j}$ must satisfy the relation
\begin{equation}\label{eq:LinearSystemWeights}
	\sum_{j=0}^{s} w_{n,j} j^{\nu} = - \sum_{j=0}^{n} \omega_{n-j} j^{\nu} + \frac{\Gamma(\nu+1)}{\Gamma(1+\nu+\alpha)} n^{\nu+\alpha}
	, \quad
	\nu \in {\cal A}_{1}\cup \bigl\{1\bigr\}
\end{equation}
and, hence, they are evaluated by solving, at each step $n$, a system of $s+1$ linear equations, with $s$ the cardinality of ${\cal A}_{1}$. 

This problem has been in--depth discussed in \cite{DiethelmFordFordWeilbeer2006}, where the authors stressed the attention on the mildly ill--conditioned nature of the coefficient matrix $V$ of (\ref{eq:LinearSystemWeights}). Nevertheless, for the methods under investigation this is not a serious issue since the size of $V$ is small, being $s=\left\lceil 1/\alpha \right\rceil + 1$ (in most cases $s$ will be 3 or 4 and the inversion of $V$ can be done analytically). The only exception is for very small values of $\alpha$ which, however, are extremely uncommon in applications.

\subsection{Starting the computation}

The initialization of (\ref{eq:FLMM}) requires the knowledge of the first $s+1$ approximations $y_{0}, y_{1}, \dots y_{s}$ of the solution. Since just $y_{0}$ is given by the problem, the remaining values must be evaluated in some other way, for instance by using some different methods.

To avoid mixing methods of different nature, we prefer to evaluate all the approximations $y_{1}, \dots y_{s}$ at the same time by setting the single implicit system
\[
	\left( \begin{array}{c}
		y_{1} \\ y_{2} \\ \vdots \\ y_{s} \
	\end{array} \right) 
	= \left( \begin{array}{c}
				T_{m-1}(t_{1}) \\ T_{m-1}(t_{2}) \\ \vdots \\ T_{m-1}(t_{s}) \
			\end{array} \right) 
			+
			h^{\alpha} 
			\left( \begin{array}{c}
				( \omega_{1} + w_{1,0} ) f_{0} \\ ( \omega_{2} + w_{2,0} ) f_{0} \\ \vdots \\ ( \omega_{s} + w_{s,0} ) f_{0} \
			\end{array} \right) \\
	+ h^{\alpha} {\cal B}
			\left( \begin{array}{c}
				f(t_{1},y_{1}) \\ f(t_{2},y_{2}) \\ \vdots \\ f(t_{s},y_{s}) \
			\end{array} \right) ,
\]
where
\[
	{\cal B} = 
			\left( \begin{array}{cccc}
				\omega_{0} I & & & \\
				\omega_{1} I & \omega_{0} I & & \\
				\vdots & \vdots & \ddots & \vdots \\
				\omega_{s-1} I & \omega_{s-2} I & \dots & \omega_{0} I \
			\end{array} \right) 
			+
			\left( \begin{array}{cccc}
				w_{1,1} I & w_{1,2} I & \dots & w_{1,s} I\\
				w_{2,1} I & w_{2,2} I & \dots & w_{2,s} I\\
				\vdots & \vdots & \ddots & \vdots \\
				w_{s,1} I & w_{s,2} I & \dots & w_{s,s} I\\
			\end{array} \right) .
\]

This nonlinear system has size $s \cdot q$, where $q$ is the size of the system of FDEs we want to solve. No particular problems arise when the size $q$ of the problem is small; otherwise some memory allocation problems could occur in presence of very small values of $\alpha$, involving a very large number $s$ of initial  approximations to compute. 

\subsection{Solution of the nonlinear systems}

All the implicit methods analyzed in Sections \ref{S:ProductIntegration} and \ref{S:FLMM} require, at each step, the  solution of an algebraic nonlinear system in the form
\[
	y_{n} = g_{n} + h^{\alpha} \bar{\omega}_{0} f(t_{n},y_{n}) ,
\]
where $g_{n}$ and $\bar{\omega}_{0}$ do not depend on $y_{n}$. To preserve the good stability properties, it is convenient to use Newton iterations 
\[
	y_{n}^{(k+1)} = y_{n}^{(k)} + \Delta^{(k)}_{n} ,
\]
where $\Delta^{(k)}_{n}$ is solution of the linear system
\[
	A \Delta^{(k)}_{n} = b
	, \quad
	A = I - h^{\alpha} \bar{\omega}_{0} J_{f}(y_{n}^{(k)}) 
	, \quad
	b = - y_{n}^{(k)} + g_{n} + h^{\alpha} \bar{\omega}_{0} f(t_{n},y_{n}^{(k)}) 
\]
and, usually, $y_{n}^{(0)} = y_{n-1}$. We refer to \cite{HairerWanner1996} for more insights on this issue and for a discussion on the stopping criteria.

We just highlight here that the convergence can depend on $h^{\alpha}$, instead of $h$ as in ODEs. Whilst this is not an issue when $1<\alpha<2$, some unexpected problems could arise whenever $\alpha$ is small. For this reason it is necessary to use a quite small step--size $h$ with low values of $\alpha$.

\section{Numerical experiments}\label{S:NumericalExperiments}

By means of numerical experiments we compare the methods discussed in the paper. We perform all the experiments in Matlab, version 7.9.0.529, on a computer equipped with the Intel Dual Core E5400 processor running at 2.70 GHz under the Windows XP operating system; all the Matlab codes have been optimized, to get the best performance, according to the discussion reported in the previous sections. %All the codes are freely available on the web page of the author.

%We have tried to optimize the Matalb codes according to the discussion reported in the previous sections. All the codes are freely available on the web page of the author.

As a first test problem we consider the linear FDE (\ref{eq:FDE_Linear}) for $\lambda=-2$. For each method, and for an increasing number $N$ of grid--points, we report the error $E(N)$ at $T=2$, with respect to a reference solution, and we provide an estimation of the order of convergence (EOC) as $\log_{2}(E(N)/E(2N))$.

\begin{table}[htb]
\footnotesize
\[
   \begin{array}{|r|cc|cc|cc|cc|cc|} \hline
& \multicolumn{2}{|c|}{\textrm{PI U}}& \multicolumn{2}{|c|}{\textrm{PI G}}& \multicolumn{2}{|c|}{\textrm{FT}}& \multicolumn{2}{|c|}{\textrm{NG}}& \multicolumn{2}{|c|}{\textrm{FBDF}}\\ 
N & \textrm{Error} & \textrm{EOC} & \textrm{Error} & \textrm{EOC} & \textrm{Error} & \textrm{EOC} & \textrm{Error} & \textrm{EOC} & \textrm{Error} & \textrm{EOC} \\ \hline
32 & 3.29(-4) & & 1.45(-4) & & 1.71(-5) & & 3.92(-5) & & 1.10(-4) & \\ 
64 & 1.15(-4) & 1.524 & 3.65(-5) & 1.987 & 5.65(-6) & 1.602 & 1.20(-5) & 1.707 & 3.16(-5) & 1.798 \\ 
128 & 4.00(-5) & 1.516 & 9.17(-6) & 1.991 & 1.74(-6) & 1.698 & 3.50(-6) & 1.780 & 8.83(-6) & 1.842 \\ 
256 & 1.40(-5) & 1.511 & 2.30(-6) & 1.993 & 5.07(-7) & 1.779 & 9.78(-7) & 1.838 & 2.40(-6) & 1.880 \\ 
512 & 4.94(-6) & 1.508 & 5.78(-7) & 1.994 & 1.41(-7) & 1.844 & 2.65(-7) & 1.883 & 6.37(-7) & 1.912 \\ 
1024 & 1.74(-6) & 1.505 & 1.45(-7) & 1.992 & 3.77(-8) & 1.903 & 6.98(-8) & 1.924 & 1.66(-7) & 1.939 \\ 
2048 & 6.14(-7) & 1.503 & 3.67(-8) & 1.987 & 9.49(-9) & 1.991 & 1.77(-8) & 1.978 & 4.25(-8) & 1.969 \\ 
\hline
   \end{array}
\]
\normalsize
\caption{Errors and EOC at $T=2.0$ for the linear FDE (\ref{eq:FDE_Linear}) with $\alpha=0.5$ and $\lambda=-2.0$} 
\label{tab:Errors_probl1_alpha050} 
\end{table}

The numerical results in Table \ref{tab:Errors_probl1_alpha050} for $\alpha=0.5$ clearly confirm theoretical findings on the order of convergence. As expected from Theorem \ref{thm:PI_Convergence}, the PI rule on uniform grids fails to obtain order 2 and converges with order $1+\alpha$; on the other hand, the PI rule on graded grids converges with full order $2$. 

As in the ODE case, FT produces the smallest error. Although also NG and PI are generalizations of the trapezoidal rule, it seems that only FT inherits this feature; by taking into account, from the analysis of Section \ref{S:LinearStability}, that FT preserves the other peculiar character of the trapezoidal rule of having the region of stability exactly equal to the region in which the solution converges towards the steady state, it is tempting to conclude that FT is the most natural generalization of the trapezoidal rule.

Unfortunately, the situation is more complicate since when $1<\alpha<2$ it is NG that provides the lowest error, as shown in Table \ref{tab:Errors_probl1_alpha150} for $\alpha=1.5$. 

\begin{table}[htb]
\footnotesize
\[
   \begin{array}{|r|cc|cc|cc|cc|cc|} \hline
& \multicolumn{2}{|c|}{\textrm{PI U}}& \multicolumn{2}{|c|}{\textrm{PI G}}& \multicolumn{2}{|c|}{\textrm{FT}}& \multicolumn{2}{|c|}{\textrm{NG}}& \multicolumn{2}{|c|}{\textrm{FBDF}}\\ 
N & \textrm{Error} & \textrm{EOC} & \textrm{Error} & \textrm{EOC} & \textrm{Error} & \textrm{EOC} & \textrm{Error} & \textrm{EOC} & \textrm{Error} & \textrm{EOC} \\ \hline
64 & 3.71(-5) & & 6.16(-5) & & 5.50(-5) & & 1.55(-5) & & 1.95(-4) & \\ 
128 & 9.31(-6) & 1.993 & 1.54(-5) & 2.000 & 1.39(-5) & 1.988 & 3.73(-6) & 2.053 & 5.22(-5) & 1.902 \\ 
256 & 2.33(-6) & 1.997 & 3.85(-6) & 2.001 & 3.48(-6) & 1.993 & 9.10(-7) & 2.035 & 1.35(-5) & 1.951 \\ 
512 & 5.82(-7) & 2.003 & 9.59(-7) & 2.004 & 8.71(-7) & 1.999 & 2.22(-7) & 2.035 & 3.43(-6) & 1.975 \\ 
1024 & 1.43(-7) & 2.024 & 2.37(-7) & 2.016 & 2.16(-7) & 2.013 & 5.25(-8) & 2.081 & 8.65(-7) & 1.989 \\ 
2048 & 3.32(-8) & 2.106 & 5.67(-8) & 2.064 & 5.15(-8) & 2.068 & 1.05(-8) & 2.325 & 2.15(-7) & 2.008 \\ 
\hline
   \end{array}
\]
\normalsize
\caption{Errors and EOC at $T=2.0$ for the linear FDE (\ref{eq:FDE_Linear}) with $\alpha=1.5$ and $\lambda=-2.0$} 
\label{tab:Errors_probl1_alpha150} 
\end{table}

\begin{figure}[htb]
	\centering
	\includegraphics[width=0.73\textwidth]{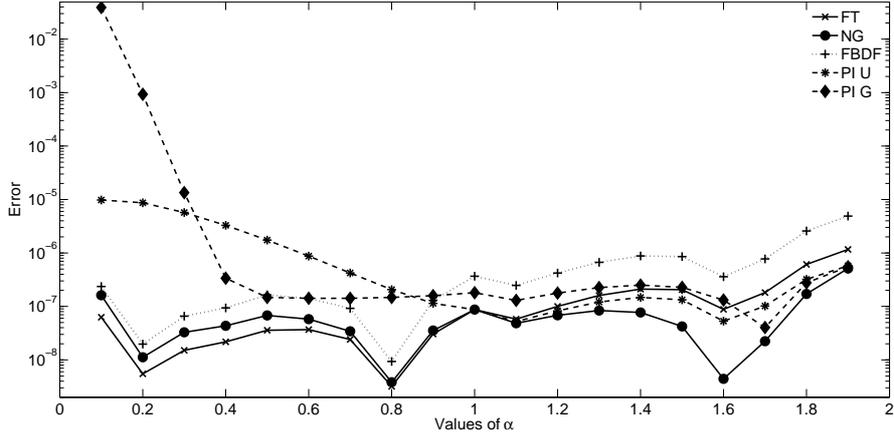} 
\caption{Errors at $T=2.0$ for the linear FDE (\ref{eq:FDE_Linear}) with $\lambda=-2.0$ and $N=1024$.}
\label{fig:Fig_ErrorALpha}
\end{figure}

This phenomenon can be more clearly observed in Figure \ref{fig:Fig_ErrorALpha}, where we show the errors as a function of $\alpha$ (we used the same values for $\lambda$ and $T$ from the previous experiment and a fixed number $N=1024$ of grid points).

%It is easily seen that whilst FT presents the lowest error for $0<\alpha<1$, it is NG to obtain this performance when $1<\alpha<2$. 
We must highlight that this is a not completely unexpected result. Stability and accuracy are indeed conflicting requirements in most of numerical methods; thus, FT shows the smallest regions of stability but the lowest error when $0<\alpha<1$, exactly as it happens with NG for $1<\alpha<2$. What we observe is a kind of roles swap between FT and NG when $\alpha$ moves from the left to the right of 1.

The very bad performance of PI on graded grids for low values of $\alpha$ is due to the fact that the grids generated by the graded exponent $r=2\alpha$ are very coarse close to right endpoint of the integration interval. Furthermore, graded grids do not add any improvement for $\alpha >1$, since the maximum expected order $2$ is already achieved by using uniform grids. We can draw the conclusion that the utility of graded grids is confined only to $\frac{1}{2}<\alpha<1$.

\begin{figure}[htb]
	\centering
	\begin{tabular}{c@{\hspace{0.8cm}}c}
		\includegraphics[width=0.40\textwidth,height=0.40\textwidth]{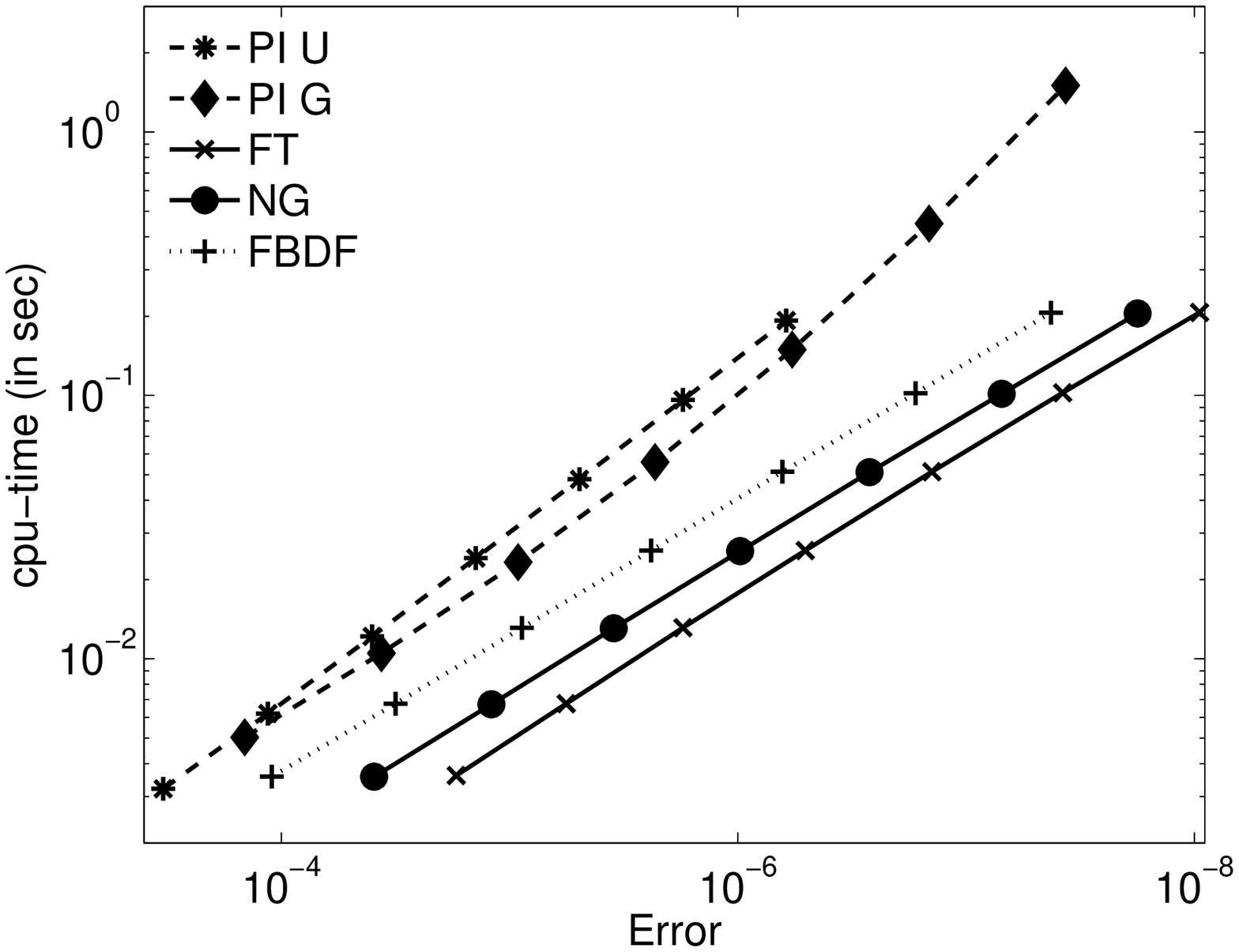} 
		&
		\includegraphics[width=0.40\textwidth,height=0.40\textwidth]{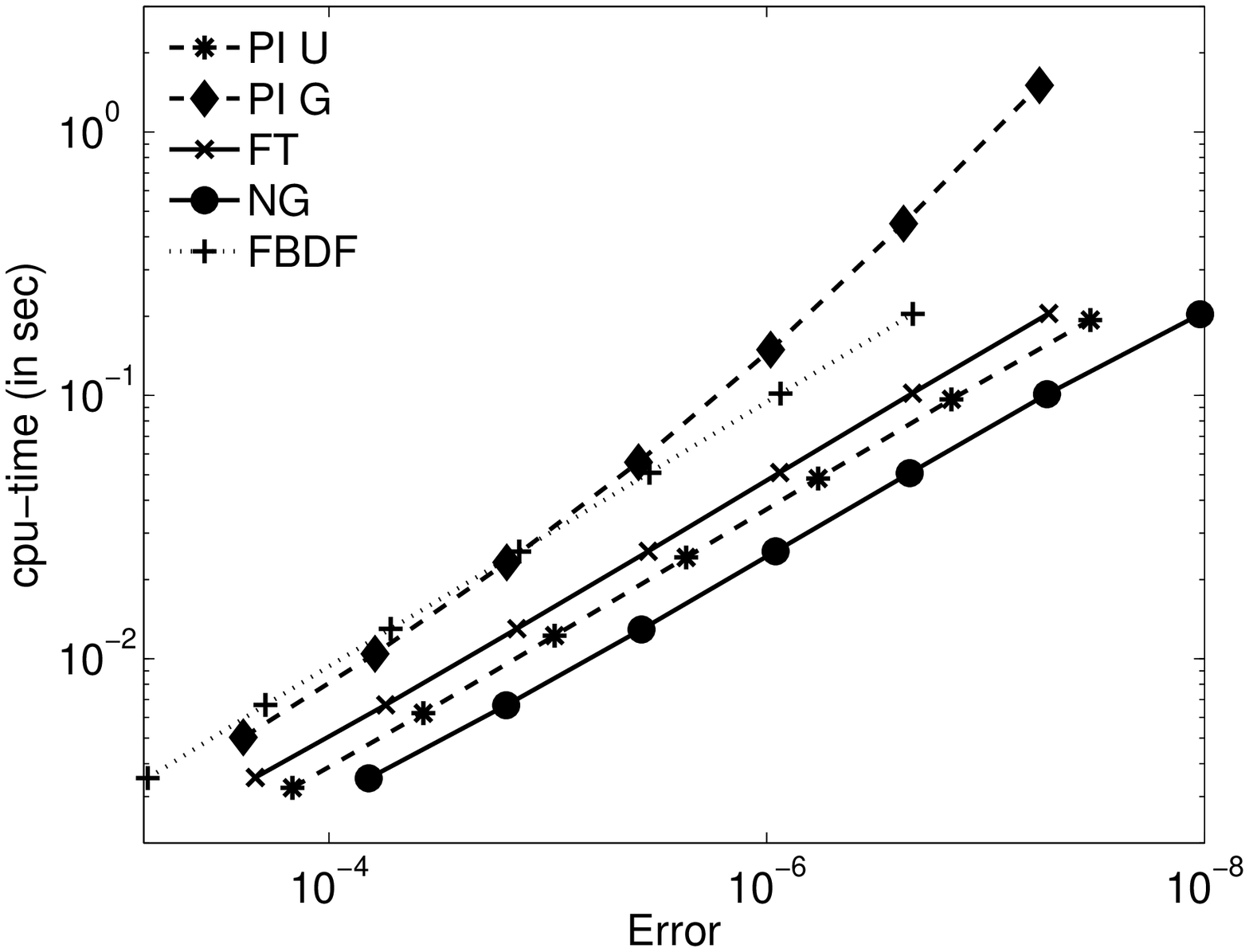} 
	\end{tabular}
\caption{Errors versus time of execution  for $\alpha=0.5$ (left) and $\alpha=1.5$ (right)}
\label{fig:Fig_ErrorTime_Probl1}
\end{figure}

It is also interesting to analyze the computational effort in relation to the accuracy. To this purpose, the plots in Figure \ref{fig:Fig_ErrorTime_Probl1}, which are related to the experiments presented in Tables \ref{tab:Errors_probl1_alpha050} and \ref{tab:Errors_probl1_alpha150}, show that FT, NG and FBDF behave in a more efficient way with respect to PI when $0<\alpha<1$; indeed, PI U pays for the reduction in the accuracy whilst PI G pays for the absence of a convolution structure which precludes the application of FFT algorithms. The minor efficiency of FBDF with respect to FT and NG is due to the major complexity in the evaluation of weights. It is worthwhile to note that the preeminence of FLMMs with respect to PI rules holds despite the need of solving linear systems for evaluating the starting weights. For $1<\alpha <2$ the absence of order reduction increases the efficiency of PI on uniform nodes but however it is below that of NG which instead has lower errors.

To observe the behavior of the methods under investigation for nonlinear problems, we consider the fractional reaction diffusion system of Brusselator type 
\begin{equation}\label{eq:Brusselator}
	\left\{ \begin{array}{l}
	{}^C D^{\alpha}_{t_{0}} x_{1}(t) = a-(\mu +1)x_{1}(t)+x_{1}(t)^{2}x_{2}(t) \\ 
	{}^C D^{\alpha}_{t_{0}} x_{2}(t) = \mu x_{1}(t)-x_{1}(t)^{2}x_{2}(t) \ 
	\end{array} \right. .
\end{equation}

\begin{figure}[htb]
	\centering
	\begin{tabular}{c@{\hspace{0.8cm}}c}
		\includegraphics[width=0.40\textwidth,height=0.40\textwidth]{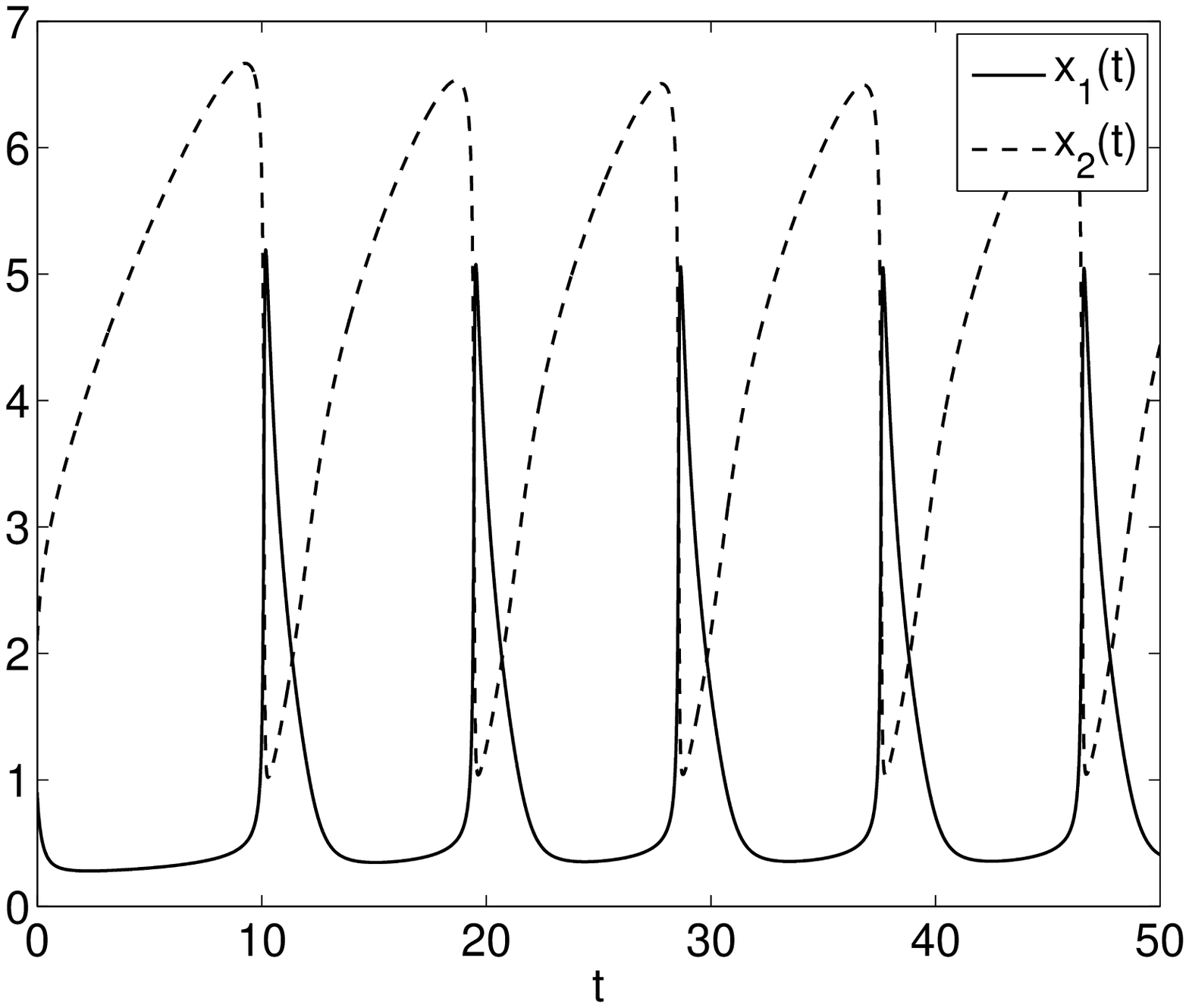} 
		&
		\includegraphics[width=0.40\textwidth,height=0.40\textwidth]{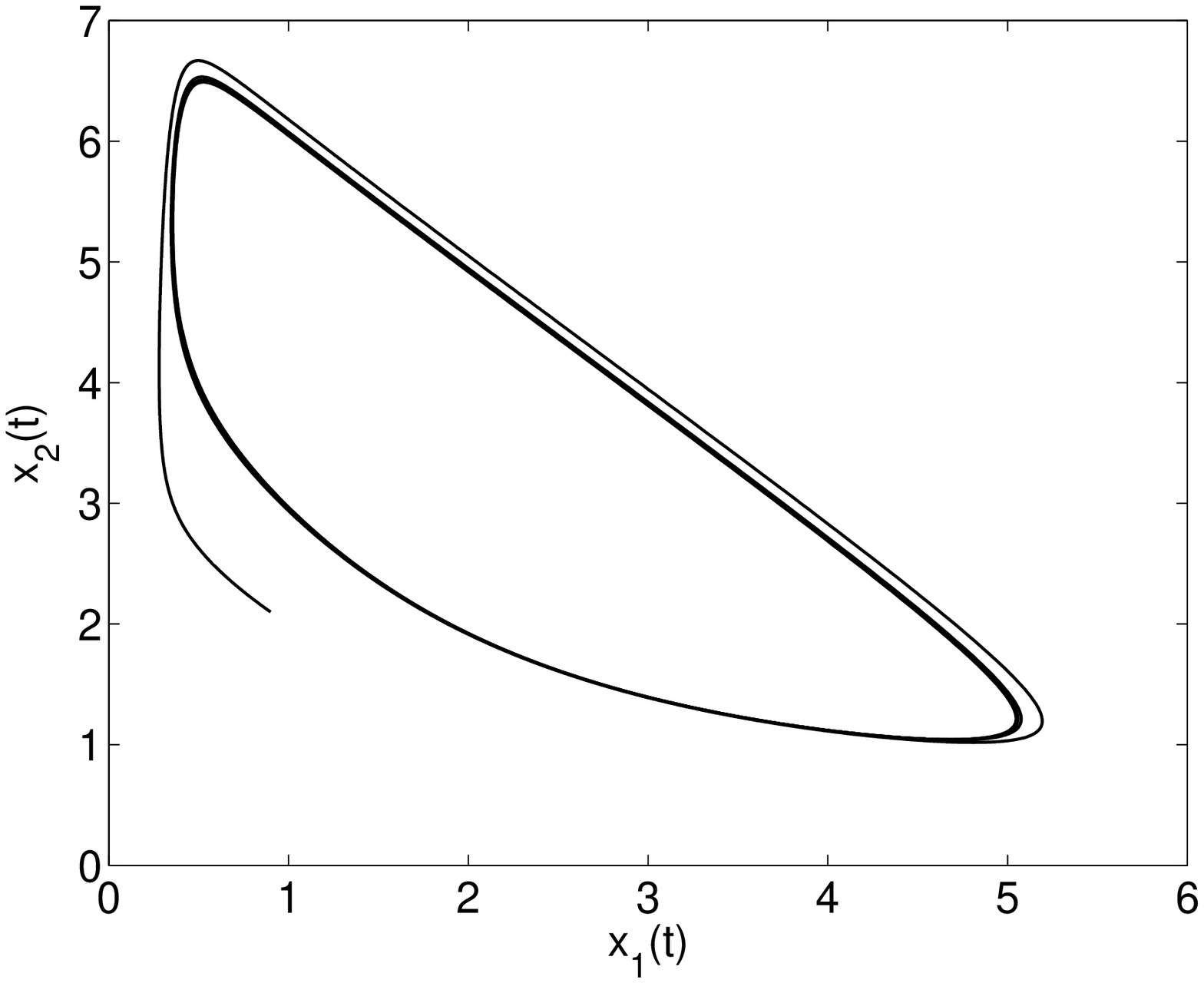} 
	\end{tabular}
\caption{Solution of the fractional Brusselator system (\ref{eq:Brusselator}) for $\alpha=0.8$ and $(a,\mu) = (1,4)$ in the $(t,x)$-plane (left) and in the phase plane (right)}
\label{fig:Fig_Soluz_Brusselator}
\end{figure}

It has been shown \cite{GafiychukDatsko2008} that there exists a marginal value $\bar{\alpha}$, depending on $a$ and $\mu$, such that system (\ref{eq:Brusselator}) approaches a stable limit cycle, as shown in Figure \ref{fig:Fig_Soluz_Brusselator} for $\alpha=0.8$ and $(a,\mu) = (1,4)$.

\begin{table}[htb]
\footnotesize
\[
   \begin{array}{|r|cc|cc|cc|cc|cc|} \hline
& \multicolumn{2}{|c|}{\textrm{PI U}}& \multicolumn{2}{|c|}{\textrm{PI G}}& \multicolumn{2}{|c|}{\textrm{FT}}& \multicolumn{2}{|c|}{\textrm{NG}}& \multicolumn{2}{|c|}{\textrm{FBDF}}\\ 
N & \textrm{Error} & \textrm{EOC} & \textrm{Error} & \textrm{EOC} & \textrm{Error} & \textrm{EOC} & \textrm{Error} & \textrm{EOC} & \textrm{Error} & \textrm{EOC} \\ \hline
400 & 3.66(-1) & & 7.33(-1) & & 2.88(-1) & & 3.67(-1) & & 5.23(-1) & \\ 
800 & 6.04(-2) & 2.601 & 2.55(-1) & 1.523 & 5.12(-2) & 2.490 & 6.32(-2) & 2.537 & 1.50(-1) & 1.800 \\ 
1600 & 1.55(-2) & 1.965 & 5.66(-2) & 2.173 & 1.28(-2) & 2.005 & 1.63(-2) & 1.956 & 4.31(-2) & 1.798 \\ 
3200 & 4.00(-3) & 1.953 & 1.35(-2) & 2.064 & 3.27(-3) & 1.964 & 4.21(-3) & 1.951 & 1.19(-2) & 1.864 \\ 
6400 & 1.01(-3) & 1.982 & 3.37(-3) & 2.006 & 8.27(-4) & 1.982 & 1.07(-3) & 1.982 & 3.12(-3) & 1.926 \\ 
12800 & 2.41(-4) & 2.071 & 8.33(-4) & 2.015 & 1.94(-4) & 2.095 & 2.55(-4) & 2.065 & 7.92(-4) & 1.977 \\ 
\hline
   \end{array}
\]
\normalsize
\caption{Errors and EOC at $T=50.0$ for the Brusselator system (\ref{eq:Brusselator}) with $\alpha=0.8$} 
\label{tab:Errors_probl3_alpha080} 
\end{table}

From Table \ref{tab:Errors_probl3_alpha080} the good results provided by PI U stand out. This is an indirect evidence that the findings of Theorem \ref{thm:PI_Convergence} can be quite pessimistic in some cases; indeed, the drop in the convergence order depends on the expansion of the vector field which in some cases could not involve real powers of low degree. However, also in this case FT presents the smallest error.

\begin{figure}[htb]
	\centering
	\begin{tabular}{c}
		\includegraphics[width=0.40\textwidth,height=0.40\textwidth]{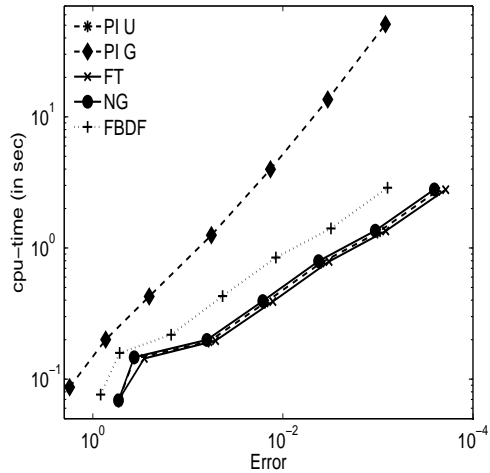} 
	\end{tabular}
\caption{Errors versus time of execution for the Brusselator system (\ref{eq:Brusselator}) with $\alpha=0.8$}
\label{fig:Fig_ErrorTime_Probl3}
\end{figure}

The relationship between errors and execution times in Figure \ref{fig:Fig_ErrorTime_Probl3} shows however that the efficiency of PI U, also in this favourable case, can be at least comparable but not higher with respect to FT and NG; it is hence quite difficult to justify the major attention that PI rules usually attract in the literature on the numerical solution of FDEs.

\section{Concluding remarks}

In this paper we have discussed some implicit second order methods for FDEs and investigated some of the main aspects related to their efficient implementation.

We have observed that stability properties change according to the value of the fractional order $\alpha$ and whenever $1<\alpha<2$ the NG and PI are no longer $A(\alpha\frac{\pi}{2})$--stable. As in the ODE case, the method obtained starting from BDF shows the largest stability regions.

Whilst FT provides the smallest error for $0<\alpha<1$, the performance of NG stands out for $1<\alpha<2$ . PI rules with graded grids appear effective just for $\alpha$ between $\frac{1}{2}$ and $1$; anyway, it seems that in general the overall efficiency of FLMMs is superior to PI, despite the major interest that PI rules usually attract.

%\bibliographystyle{elsart-num-sort}
%\bibliographystyle{acm}
%\bibliography{FractionalTrapezoidal}

\begin{thebibliography}{10}

\bibitem{Brunner1985}
{\sc Brunner, H.}
\newblock The numerical solution of weakly singular {V}olterra integral
  equations by collocation on graded meshes.
\newblock {\em Math. Comp. 45}, 172 (1985), 417--437.

\bibitem{Brunner2004}
{\sc Brunner, H.}
\newblock {\em Collocation methods for {V}olterra integral and related
  functional differential equations}.
\newblock Cambridge University Press, Cambridge, 2004.

\bibitem{BuenoOrovioKayGrauRodriguezBurrage2013}
{\sc Bueno-Orovio, A., Kay, D., Grau, V., Rodriguez, B., and Burrage, K.}
\newblock Fractional diffusion models of cardiac electrical propagation : role
  of structural heterogeneity in dispersion of repolarization.
\newblock Tech. Rep. OCCAM 13/35, Oxford Centre for Collaborative Applied
  Mathematics, Oxford (UK), 2013.

\bibitem{BurrageHaleKay2010}
{\sc Burrage, K., Hale, N., and Kay, D.}
\newblock An efficient implicit {FEM} scheme for fractional-in-space
  reaction-diffusion equations.
\newblock {\em SIAM J. Sci. Comput. 34}, 4 (2012), A2145--A2172.

\bibitem{CafagnaGrassi2012}
{\sc Cafagna, D., and Grassi, G.}
\newblock Observer-based projective synchronization of fractional systems via a
  scalar signal: Application to hyperchaotic {R}\"ossler systems.
\newblock {\em Nonlinear Dynamics 68}, 1-2 (2012), 117--128.

\bibitem{CameronMcKee1985}
{\sc Cameron, R.~F., and McKee, S.}
\newblock The analysis of product integration methods for {A}bel's equation
  using discrete fractional differentiation.
\newblock {\em IMA J. Numer. Anal. 5}, 1 (1985), 339--353.

\bibitem{CaponettoMaionePisanoRapaiUsai2013}
{\sc Caponetto, R., Maione, G., Pisano, A., Rapai{\'c}, M. M.~R., and Usai, E.}
\newblock Analysis and shaping of the self-sustained oscillations in relay
  controlled fractional-order systems.
\newblock {\em Fractional Calculus and Applied Analysis 16}, 1 (2013), 93--108.

\bibitem{ConteDelPrete2006}
{\sc Conte, D., and Prete, I.~D.}
\newblock Fast collocation methods for {V}olterra integral equations of
  convolution type.
\newblock {\em J. Comput. Appl. Math. 196}, 2 (2006), 652--663.

\bibitem{HoogWeiss1974}
{\sc de~Hoog, F., and Weiss, R.}
\newblock High order methods for a class of {V}olterra integral equations with
  weakly singular kernels.
\newblock {\em SIAM J. Numer. Anal. 11\/} (1974), 1166--1180.

\bibitem{DieciLopez2009}
{\sc Dieci, L., and Lopez, L.}
\newblock Sliding motion in {F}ilippov differential systems: theoretical
  results and a computational approach.
\newblock {\em SIAM J. Numer. Anal. 47}, 3 (2009), 2023--2051.

\bibitem{DieciLopez2011}
{\sc Dieci, L., and Lopez, L.}
\newblock Sliding motion on discontinuity surfaces of high co-dimension. {A}
  construction for selecting a {F}ilippov vector field.
\newblock {\em Numer. Math. 117}, 4 (2011), 779--811.

\bibitem{DieciLopez2012}
{\sc Dieci, L., and Lopez, L.}
\newblock A survey of numerical methods for {IVP}s of {ODE}s with discontinuous
  right-hand side.
\newblock {\em J. Comput. Appl. Math. 236}, 16 (2012), 3967--3991.

\bibitem{Diethelm2008}
{\sc Diethelm, K.}
\newblock An investigation of some nonclassical methods for the numerical
  approximation of {C}aputo-type fractional derivatives.
\newblock {\em Numer. Algorithms 47}, 4 (2008), 361--390.

\bibitem{Diethelm2010}
{\sc Diethelm, K.}
\newblock {\em The analysis of fractional differential equations}, vol.~2004 of
  {\em Lecture Notes in Mathematics}.
\newblock Springer-Verlag, Berlin, 2010.

\bibitem{Diethelm2011}
{\sc Diethelm, K.}
\newblock An efficient parallel algorithm for the numerical solution of
  fractional differential equations.
\newblock {\em Fract. Calc. Appl. Anal. 14}, 3 (2011), 475--490.

\bibitem{DiethelmFordFordWeilbeer2006}
{\sc Diethelm, K., Ford, J.~M., Ford, N.~J., and Weilbeer, M.}
\newblock Pitfalls in fast numerical solvers for fractional differential
  equations.
\newblock {\em J. Comput. Appl. Math. 186}, 2 (2006), 482--503.

\bibitem{DiethelmFordFreedLuchko2005}
{\sc Diethelm, K., Ford, N.~J., Freed, A.~D., and Luchko, Y.}
\newblock Algorithms for the fractional calculus: a selection of numerical
  methods.
\newblock {\em Comput. Methods Appl. Mech. Engrg. 194}, 6-8 (2005), 743--773.

\bibitem{DiethelmFreed1999}
{\sc Diethelm, K., and Freed, A.~D.}
\newblock The {F}rac{P}{E}{C}{E} subroutine for the numerical solution of
  differential equations of fractional order.
\newblock In {\em Forschung und wissenschaftliches Rechnen 1998}, S.Heinzel and
  T.Plesser, Eds. 1999, pp.~57--71.

\bibitem{Dixon1985}
{\sc Dixon, J.}
\newblock On the order of the error in discretization methods for weakly
  singular second kind {V}olterra integral equations with nonsmooth solutions.
\newblock {\em BIT 25}, 4 (1985), 624--634.

\bibitem{GafiychukDatsko2008}
{\sc Gafiychuk, V.~V., and Datsko, B.}
\newblock Stability analysis and limit cycle in fractional system with
  {B}russelator nonlinearities.
\newblock {\em Physics Letters A 372}, 29 (2008), 4902 -- 4904.

\bibitem{GaleoneGarrappa2008}
{\sc Galeone, L., and Garrappa, R.}
\newblock Fractional {A}dams-{M}oulton methods.
\newblock {\em Math. Comput. Simulation 79}, 4 (2008), 1358--1367.

\bibitem{Garra2011}
{\sc Garra, R.}
\newblock Fractional-calculus model for temperature and pressure waves in
  fluid-saturated porous rocks.
\newblock {\em Phys. Rev. E 84\/} (Sep 2011), 036605.

\bibitem{Garrappa2010_IJCM}
{\sc Garrappa, R.}
\newblock On linear stability of predictor--corrector algorithms for fractional
  differential equations.
\newblock {\em Int. J. Comput. Math. 87}, 10 (2010), 2281--2290.

\bibitem{Garrappa2012_MCS}
{\sc Garrappa, R.}
\newblock On some generalizations of the implicit {E}uler method for
  discontinuous fractional differential equations.
\newblock {\em Math. Comput. Simulat.\/} (2012).

\bibitem{GarrappaPopolizio2011_JCAM}
{\sc Garrappa, R., and Popolizio, M.}
\newblock On accurate product integration rules for linear fractional
  differential equations.
\newblock {\em J. Comput. Appl. Math. 235}, 5 (2011), 1085--1097.

\bibitem{GarrappaPopolizio2013}
{\sc Garrappa, R., and Popolizio, M.}
\newblock Evaluation of generalized {M}ittag--{L}effler functions on the real
  line.
\newblock {\em Adv. Comput. Math. 39}, 1 (2013), 205--225.

\bibitem{Gorenflo1997}
{\sc Gorenflo, R.}
\newblock Fractional calculus: some numerical methods.
\newblock In {\em Fractals and Fractional Calculus in Continuum Mechanics},
  A.~Carpinteri and F.~Mainardi, Eds., vol.~278 of {\em CISM Courses and
  Lectures}. Springer Verlag, Wien and New York, 1997, pp.~277--290.

\bibitem{HairerLubichSchlichte1985}
{\sc Hairer, E., Lubich, C., and Schlichte, M.}
\newblock Fast numerical solution of nonlinear {V}olterra convolution
  equations.
\newblock {\em SIAM J. Sci. Statist. Comput. 6}, 3 (1985), 532--541.

\bibitem{HairerWanner1996}
{\sc Hairer, E., and Wanner, G.}
\newblock {\em Solving ordinary differential equations. {II}}, second~ed.
\newblock Springer-Verlag, Berlin, 1996.

\bibitem{Henrici1974}
{\sc Henrici, P.}
\newblock {\em Applied and computational complex analysis}, vol.~1.
\newblock John Wiley \& Sons, New York, 1974.

\bibitem{Henrici1979}
{\sc Henrici, P.}
\newblock Fast {F}ourier methods in computational complex analysis.
\newblock {\em SIAM Rev. 21}, 4 (1979), 481--527.

\bibitem{KademLuchkoBaleanu2010}
{\sc Kadem, A., Luchko, Y., and Baleanu, D.}
\newblock Spectral method for solution of the fractional transport equation.
\newblock {\em Rep. Math. Phys. 66}, 1 (2010), 103--115.

\bibitem{KilbasSrivastavaTrujillo2006}
{\sc Kilbas, A.~A., Srivastava, H.~M., and Trujillo, J.~J.}
\newblock {\em Theory and applications of fractional differential equations},
  vol.~204 of {\em North-Holland Mathematics Studies}.
\newblock Elsevier Science B.V., Amsterdam, 2006.

\bibitem{Lambert1991}
{\sc Lambert, J.~D.}
\newblock {\em Numerical methods for ordinary differential systems}.
\newblock John Wiley \& Sons Ltd., Chichester, 1991.

\bibitem{LiZeng2013}
{\sc Li, C., and Zeng, F.}
\newblock The finite difference methods for fractional ordinary differential
  equations.
\newblock {\em Numer. Funct. Anal. Optim. 34}, 2 (2013), 149--179.

\bibitem{Lubich1983}
{\sc Lubich, C.}
\newblock Runge-{K}utta theory for {V}olterra and {A}bel integral equations of
  the second kind.
\newblock {\em Math. Comp. 41}, 163 (1983), 87--102.

\bibitem{Lubich1985}
{\sc Lubich, C.}
\newblock Fractional linear multistep methods for {A}bel--{V}olterra integral
  equations of the second kind.
\newblock {\em Math. Comp. 45}, 172 (1985), 463--469.

\bibitem{Lubich1986}
{\sc Lubich, C.}
\newblock Discretized fractional calculus.
\newblock {\em SIAM J. Math. Anal. 17}, 3 (1986), 704--719.

\bibitem{Lubich1986_IMA}
{\sc Lubich, C.}
\newblock A stability analysis of convolution quadratures for {A}bel-{V}olterra
  integral equations.
\newblock {\em IMA J. Numer. Anal. 6}, 1 (1986), 87--101.

\bibitem{MachadoKiryakovaMainardi2011}
{\sc Machado, J.~T., Kiryakova, V., and Mainardi, F.}
\newblock Recent history of fractional calculus.
\newblock {\em Commun. Nonlinear Sci. Numer. Simul. 16}, 3 (2011), 1140--1153.

\bibitem{Magin2010}
{\sc Magin, R.~L.}
\newblock Fractional calculus models of complex dynamics in biological tissues.
\newblock {\em Computers \& Mathematics with Applications 59}, 5 (2010), 1586
  -- 1593.

\bibitem{Mainardi2010}
{\sc Mainardi, F.}
\newblock {\em Fractional calculus and waves in linear viscoelasticity}.
\newblock Imperial College Press, London, 2010.

\bibitem{Matignon1998}
{\sc Matignon, D.}
\newblock Stability properties for generalized fractional differential systems.
\newblock In {\em Syst\`emes diff\'erentiels fractionnaires ({P}aris, 1998)},
  vol.~5 of {\em ESAIM Proc.} Soc. Math. Appl. Indust., Paris, 1998,
  pp.~145--158.

\bibitem{MeerschaertTadjeran2004}
{\sc Meerschaert, M.~M., and Tadjeran, C.}
\newblock Finite difference approximations for fractional advection-dispersion
  flow equations.
\newblock {\em J. Comput. Appl. Math. 172}, 1 (2004), 65--77.

\bibitem{MeerschaertTadjeran2006}
{\sc Meerschaert, M.~M., and Tadjeran, C.}
\newblock Finite difference approximations for two-sided space-fractional
  partial differential equations.
\newblock {\em Appl. Numer. Math. 56}, 1 (2006), 80--90.

\bibitem{Moret2013}
{\sc Moret, I.}
\newblock A note on {K}rylov methods for fractional evolution problems.
\newblock {\em Numer. Func. Anal. Opt. 34}, 5 (2013), 539--556.

\bibitem{MoretNovati2011}
{\sc Moret, I., and Novati, P.}
\newblock On the convergence of {K}rylov subspace methods for matrix
  {M}ittag-{L}effler functions.
\newblock {\em SIAM J. Numer. Anal. 49}, 5 (2011), 2144--2164.

\bibitem{MoretPopolizio2012}
{\sc Moret, I., and Popolizio, M.}
\newblock The restarted shift-and-invert {K}rylov method for matrix functions.
\newblock {\em Numerical Linear Algebra with Applications\/} (2012).

\bibitem{Podlubny1999}
{\sc Podlubny, I.}
\newblock {\em Fractional differential equations}, vol.~198 of {\em Mathematics
  in Science and Engineering}.
\newblock Academic Press Inc., San Diego, CA, 1999.

\bibitem{MachadoStefanescuTintareanuBaleanu2013}
{\sc Tenreiro~Machado, J., Stefanescu, P., Tintareanu, O., and Baleanu, D.}
\newblock Fractional calculus analysis of the cosmic microwave background.
\newblock {\em Romanian Reports in Physics 65}, 1 (2013), 316--323.

\bibitem{Wolkenfelt1979}
{\sc Wolkenfelt, P. H.~M.}
\newblock Linear multistep methods and the construction of quadrature formulae
  for volterra integral and integro-differential equations.
\newblock Tech. Rep. NW 76/79, Mathematisch Centrum,, Amsterdam (Netherlands),
  1979.

\bibitem{Young1954}
{\sc Young, A.}
\newblock Approximate product-integration.
\newblock {\em Proc. Roy. Soc. London Ser. A. 224\/} (1954), 552--561.

\bibitem{YusteAcedo2005}
{\sc Yuste, S.~B., and Acedo, L.}
\newblock An explicit finite difference method and a new von {N}eumann-type
  stability analysis for fractional diffusion equations.
\newblock {\em SIAM J. Numer. Anal. 42}, 5 (2005), 1862--1874.

\bibitem{YusteQuintanaMurillo2012}
{\sc Yuste, S.~B., and Quintana-Murillo, J.}
\newblock A finite difference method with non-uniform timesteps for fractional
  diffusion equations.
\newblock {\em Comput. Phys. Commun. 183}, 12 (2012), 2594--2600.

\end{thebibliography}

\end{document}